\documentclass[a4paper,american,11pt, reqno]{amsart}
\pdfoutput=1

\usepackage[latin1]{inputenc}

\usepackage{amsmath}
\usepackage{amssymb}
\usepackage{hyperref}
\usepackage{url}
\usepackage{babel}
\usepackage{pgf,tikz,pgfplots}
\usepackage{graphicx}

\usepackage{tikz, float} \usetikzlibrary {positioning}
\usetikzlibrary{shapes.misc}
\usetikzlibrary{arrows}

\usetikzlibrary{calc}
\usepackage{amsfonts}
\input xy
\xyoption{all}

\newcommand\xqed[1]{%
  \leavevmode\unskip\penalty9999 \hbox{}\nobreak\hfill
  \quad\hbox{#1}}
\newcommand\examp{\xqed{$\triangle$}}

\newcommand{\Z}{{\mathbb Z}}

\newcommand{\PP}{{\mathbb P}}

\newcommand{\C}{{\mathbb C}}

\newcommand{\R}{{\mathbb R}}
\newcommand{\QQ}{{\mathbb Q}}

\newcommand{\T}{{\mathbb T}}
\newcommand{\TP}{{\mathbb{T}P}}

\newcommand{\A}{{\mathcal A}}

\newcommand{\sed}{{\text{sed}}}

\newcommand{\F}{{\mathcal F}}

\newcommand{\G}{{\mathcal G}}

\newcommand{\I}{{\mathcal I}}

\renewcommand{\S}{{\mathcal S}}

\newcommand{\Aff}{\text{Aff}}

\newcommand{\Ker}{\text{Ker}}
\newcommand{\va}{\varepsilon}

\newcommand{\K}{{\mathcal K}}

\DeclareMathOperator{\rank}{\text{rank}}

\DeclareMathOperator{\relint}{\text{int}}
\DeclareMathOperator{\Gr}{\text{Gr}}

\newtheorem{thm}{Theorem}[section]
\newtheorem{defi}[thm]{Definition}

\newtheorem{definition}[thm]{Definition}
\newtheorem{prop}[thm]{Proposition}
\newtheorem{proposition}[thm]{Proposition}
\newtheorem{lemma}[thm]{Lemma}
\newtheorem{cor}[thm]{Corollary}
\newtheorem{remark}[thm]{Remark}
\newtheorem{corollary}[thm]{Corollary}
\newtheorem{conj}[thm]{Conjecture}
          {\theoremstyle{definition}
}
          {\theoremstyle{definition}
\newtheorem{exa}[thm]{Example}
\newtheorem{example}[thm]{Example}
}

\newtheorem{ques}[thm]{Question}
 \numberwithin{equation}{section}

\tikzset{%
  add/.style args={#1 and #2}{to path={%
 ($(\tikztostart)!-#1!(\tikztotarget)$)--($(\tikztotarget)!-#2!(\tikztostart)$)%
  \tikztonodes}}
} 

\newcommand{\comment}[1]{}

\begin{document}
\title[{Betti numbers of real hypersurfaces near the tropical limit}]{{Bounding the Betti numbers of real hypersurfaces near the tropical limit}}

\author[Arthur Renaudineau]{Arthur Renaudineau}
\address{Arthur Renaudineau, Univ. Lille, CNRS, UMR 8524 -
Laboratoire Paul Painlev\'e, F-59000 Lille, France.}
\email{arthur.renaudineau@univ-lille.fr}
\author[Kristin Shaw]{Kristin Shaw}
\address{Kristin Shaw, 
University of Oslo, Oslo, Norway.}
\email{krisshaw@math.uio.no}


\maketitle
\begin{abstract}
We prove a bound conjectured by Itenberg on the Betti numbers of real algebraic hypersurfaces near non-singular tropical limits. 
These bounds are given in terms of the Hodge numbers of the complexification. To prove the conjecture we introduce a real variant of tropical homology and define a filtration on the corresponding chain complex inspired by Kalinin's filtration. The spectral sequence  associated to this filtration  converges to the homology groups of the real algebraic variety and we show that the terms of the  first page are tropical homology groups with $\Z_2$-coefficients. The dimensions of these homology groups correspond to the Hodge numbers of complex projective hypersurfaces by \cite{IKMZ} and \cite{ARS}. 
The bounds  
on the Betti numbers of the real part follow, as well as a criterion to obtain a maximal variety.  We also generalise a known formula relating the signature of the complex hypersurface and the Euler characteristic of the real algebraic hypersurface, as well as Haas' combinatorial criterion for the maximality of plane curves near the tropical limit.

\end{abstract}

\section{Introduction}

A \textit{real hypersurface} $V \subset \PP^{n+1}$ of degree $d$    is a hypersurface defined by a real homogeneous polynomial $f(z_0, \dots, z_{n+1}) \in \R[z_0, \dots, z_{n+1}]$ of degree $d$. 
We let $\R V$ denote the set of   real points of  $V$ and $\C V$ denote the set of its complex points. 
The following fundamental question in real algebraic geometry can be traced back beyond Hilbert's sixteenth problem, see \cite{Hilbert}, \cite{Wilson}, \cite{DegKhar} for a survey.
\begin{ques}\label{quesintro}
For any $0\leq q\leq n$, what is the maximal possible value of the $q$-th Betti number  $$b_q(\R V) := \dim H_q(\R V; \Z_2)$$ 
 among degree $d$ non-singular real algebraic hypersurfaces $V$ in $\PP^{n+1}$?
\label{question}
\end{ques}

In 1876,  Harnack \cite{Harnack} proved for  non-singular  real  plane curves the optimal  bound $b_0(\R V) \leq g(\C V) +1$, where $g(\C V)$ denotes the genus of the  complex curve. 
Beyond the case of plane curves, no optimal bounds are known in general on the individual Betti numbers of real algebraic varieties. For example, 
in the case of non-singular  real algebraic  surfaces in $\mathbb{P}^3$, the  maximal values of the individual Betti numbers are unknown beyond degree  $5$. 
It is  known that the maximal number of connected components of a non-singular real algebraic quintic surface  is either $23, 24,$ or $25$ and  the maximal value of  the first Betti number  is either $45$ or $47$, see \cite{Orevkov2001} and \cite{ItKhar}.

In relation to higher Betti numbers, in 1980 Viro formulated the following conjecture for all  real projective surfaces. 
\begin{conj}[Viro]
If  $V$ is a  non-singular real projective surface such that $\C V$ is simply connected, then
$$
b_1(\R V)\leq h^{1,1}(\C V),
$$
where  $h^{1,1}(\C V)$ denotes the  $(1,1)$-th Hodge number of $\C V$.
\end{conj}
In general, we will denote by $h^{p,q}(\C V)$ the  $(p,q)$-th Hodge number of $\C V$.

When $V$ is the double covering of $\mathbb{P}^2$ ramified along a curve of  even degree, this conjecture is a reformulation of Ragsdale's conjecture \cite{Rags},  \cite{Viro80}. The first counterexample to her conjecture was constructed by  Itenberg  \cite{Itenberg93}.  This paved the way to various counterexamples to Viro's conjecture and to constructions of real algebraic surfaces with many connected components, for example those  in \cite{Itenberg1997}, \cite{Bihan99}, and \cite{Brugalle2006}.
 
It is still not known whether Viro's conjecture is true for surfaces which are maximal in the sense of the Smith-Thom inequality (\ref{ST}) .

There are two main directions in Question \ref{question}. The first is to prohibit topologies of a real algebraic variety, as is the case for Harnack's bound.  The second direction is to provide constructions of  real algebraic varieties with  given topology.
Viro's patchworking method provided a breakthrough in the second direction \cite{Viro84}. This technique continues to be the most powerful tool to construct real algebraic varieties in toric varieties with determined topology. Here we will restrict our attention to Viro's primitive combinatorial patchworking. 
The following  was conjectured by Itenberg around 2005, and later appeared in  \cite{ItenbergSimons}.

\begin{conj} \cite[Conjecture 2.5]{ItenbergSimons} \label{conj:main}
Let $ V$ be a real hypersurface in $\PP^{n+1}$ obtained by a primitive patchworking. Then for any integer $q= 0, \dots , n$, 
$$b_q(\R V) \leq 
\begin{cases}
h^{q, q}(\C V)   \text{ for }   q = n/2,  \\
h^{q,n-q}(\C V) + 1  \text{ otherwise}.
\end{cases} $$  
\end{conj}

In the case of real algebraic surfaces in $\mathbb{P}^3$ arising from primitive  patchworkings  the above bounds were already proven by Itenberg  \cite{Itenberg1997}, 
and  are explicitly, 
$$
b_0(\R V)\leq \binom{d-1}{3}+ 1 \qquad \text{ and} \qquad b_1(\R V) \leq \dfrac{2d^3-6d^2+7d}{3}.
$$ 
For example, real algebraic surfaces of degree $5$ arising from a primitive patchworking satisfy
$b_0(\R V)\leq 5$ and $ b_1(\R V) \leq 45.$ 
Furthermore, asymptotic analogues of the bounds in Conjecture \ref{conj:main} were proved by Itenberg and Viro in \cite{IV}.

Viro's method for patchworking applies not only to real hypersurfaces in projective space but also to hypersurfaces in more general toric varieties, see for example \cite{JJR}. Real algebraic hypersurfaces arising from primitive patchworking were later  interpreted  by Viro \cite{Viro2000} as real algebraic hypersurfaces near  non-singular tropical limits, see Definition \ref{def:neartroplim} and \cite[Section 5.3]{BIMS}. Here we will use this contemporary point of view on Viro's method and relate it to Viro's original formulation in Remark \ref{Viro}. A hypersurface near the non-singular tropical limit will always be (partially) compactified in a toric variety whose fan is a subfan of the Newton polytope of the hypersurface.  
We say a Newton polytope $\Delta$ is non-singular if the associated toric variety $Y_{\Delta}$ is non-singular.
In this paper we establish the following theorem for real algebraic hypersurfaces in compact non-singular toric varieties near a non-singular tropical limit.

\begin{thm}\label{thm:main}
Let $ V$ be a compact real algebraic hypersurface with non-singular Newton polytope and near a non-singular tropical limit.  
Then for any integer $q= 0, \dots , n$,
$$b_q(\R V) \leq 
\begin{cases}
h^{q, q}(\C V)   \text{ for }   q = n/2,  \\
h^{q,n-q}(\C V) +  h^{q,q}(\C V)  \text{ otherwise}.
\end{cases} $$  
\end{thm}

When the toric variety is projective space, Conjecture \ref{conj:main} follows directly from the following theorem and the fact that $h^{q,q}(\C V) = 1 $ for $2q \neq n$  by the Lefschetz Hyperplane Section Theorem.

\subsection{A guide to the proof of Theorem \ref{thm:main}}
To  prove  Theorem \ref{thm:main}, we use  a  tropical description of primitive patchworking in terms of real phase structures, which we present in  Section \ref{section:realtropicalhypersurfaces}. We recall the  relation to the standard version of  primitive patchworking in Remark \ref{Viro}. For an $n$-dimensional  non-singular real tropical hypersurface $X$ in a tropical toric variety $Y$ the notion of a real phase structure $\mathcal{E}$ is described in  Definition  \ref{realphase}. 
In  Section \ref{sec:signcosheaf}  we describe  a cellular cosheaf on $X$ called the sign cosheaf  $\mathcal{S}_{\mathcal{E}}$. 
A  cellular  cosheaf $\G$ on a  tropical hypersurface $X$ consists of a vector space $\G(\sigma)$ for each face $\sigma $ of $X$ together with linear maps 
$i_{\sigma \tau} \colon \G(\sigma) \to \G(\tau)$ for each inclusion of faces $\tau \subset \sigma$. These linear maps must satisfy commutativity conditions for all face relations $\rho \subset \tau_1, \tau_2 \subset \sigma$.
The cellular chain complex $C_{\bullet}(X ; \S_{\mathcal{E}})$ with coefficients in $\S_{\mathcal{E}}$ together with its homology  groups $H_{\bullet}(X;\S_{\mathcal{E}})$ are defined in Definition \ref{def:realhomo}. Note that all the cosheaves in this paper are considered over $\Z_2$ otherwise it is cleary stated. In Proposition \ref{prop:realcosheaf}, we prove that the homology groups  of the sign cosheaf are  isomorphic to the homology  groups of the real part of a real  algebraic hypersurface  near the tropical limit.
The next step of the proof is to construct a  filtration of the chain complex with coefficients in  $\S_{\mathcal{E}}$, 
\begin{equation}\label{filtration}
0\subset C_{\bullet}(X; \K_n)\subset\cdots\subset C_{\bullet}(X;\K_p)\subset\cdots\subset C_{\bullet}(X;\mathcal{S}_{\mathcal{E}}), 
\end{equation}
 where the $\K_p$'s are the collection of cellular cosheaves on $X$ from Definition \ref{def:Kp}. 
 
  The tropical homology groups,   as introduced by Itenberg, Katzarkov, Mikhakin and Zharkov \cite{IKMZ},  are also homology groups of  cosheaves on tropical varieties but with rational coefficients. Here we use a $\Z_2$-variant of  this homology theory and denote the cosheaves by $\F_p$. 
For every $p$ and each face $\tau$ of $X$,  we define linear maps   $bv_p\colon \K_p(\tau) \to\F_p(\tau)$ in Definition \ref{def:bvp}, which come from the augmented filtration of a group algebra which was highlighted by Quillen \cite{Quillen}.
It follows from Lemma \ref{lemma:isoK_p},  that these linear maps are surjective and  satisfy  $\Ker(bv_p \colon K_p(\tau) \to\F_p(\tau) ) = \K_{p+1}(\tau)$. 
 Proposition \ref{prop:exactcosheaf} shows that these linear maps commute with the cosheaf maps and thus induce morphisms of chain complexes and produce the  filtration in (\ref{filtration}). 
Then we consider the spectral sequence associated to this filtration which we denote by $(E^{\bullet}_{\bullet, \bullet}, \partial_{\bullet}).$ This spectral sequence degenerates  since  it arises from a filtration of a chain complex consisting of finite dimensional chain groups. Therefore we obtain,
$$ 
\dim H_q(X;\mathcal{S}_{\mathcal{E}}) =  \sum_{p = 0}^n \dim E^{\infty}_{q, p}.$$

By Corollary  \ref{cor:relativenadtrop} we get
$$E^1_{q,p}  = H_q(C_{\bullet}(X,\K_p)/C_{\bullet}(X,\K_{p+1})) \cong H_q(X; \F_p).$$
This establishes the following theorem which bounds the Betti numbers of real algebraic hypersurfaces near the tropical limit. 
The theorem holds not just for a compact  hypersurface  in the toric variety of its Newton polytope, but also for a non-compact hypersurface obtained by removing the intersection with any of the torus orbits of that toric variety. This includes for example, hypersurfaces near the tropical limit contained in the torus or affine space. In the statement below, the notation $b_q^{BM}$ denotes the $q$-th Betti number of the  Borel-More homology group  \cite{BorelMoore} and $H^{BM}_{q}(X; \F_p)$ denotes the Borel-Moore variant of tropical homology \cite{ARS}. 
For $\Delta$ a lattice polytope, denote by $Y_{\Delta}^o$ a partial compactification of the torus corresponding to a subfan of the dual fan of $\Delta$.
\begin{thm}\label{thm:toricvar}
Let $V$ be a real algebraic hypersurface with  Newton polytope $\Delta$ in the non-singular toric variety  $Y^o_{\Delta}$  and near  the non-singular tropical limit $X$,   then 
for all $q$  we have,
 $$b_q(\R V) \leq \sum_{p = 0}^n \dim H_{q}(X; \F_p),$$
and
 $$b_q^{BM}(\R V) \leq \sum_{p = 0}^n \dim H^{BM}_{q}(X; \F_p).$$
\end{thm}

The main ingredient for proving Theorem \ref{thm:toricvar} from Theorem  \ref{thm:main}   is to relate the dimensions of tropical homology groups with $\Z_2$-coefficients and the Hodge numbers of complex hypersurfaces.  

\begin{thm} \label{thm:hodgemod2}
Let $X$ be a non-singular compact tropical hypersurface with non-singular Newton polytope $\Delta$.
 Let $V$ be a non-singular complex hypersurface in  the  non-singular complex toric variety $Y_{\Delta}$ also with Newton polytope $\Delta$.  Then for all $p$ and $q$ we have $$\dim H^{p, q}(\C V) = \rank H_q(X;  \mathcal{F}^{\Z}_p).$$
\end{thm}

\begin{proof}
By \cite[Corollary 1.4]{ARS}, the integral tropical homology groups of the hypersurface $X$ are torsion free. Therefore, we have $$\rank H_q(X; \F_p^{\Z}) = \dim H_q(X; \F_p) = \dim H_q(X; \F^{\QQ}_p), $$
where $H_q(X; \F_p^{\Z})$ and  $H_q(X; \F_p^{\QQ})$ are the   tropical homology group with $\Z$ and $\QQ$-coefficients, respectively. 
 By    \cite[Corollary 2]{IKMZ} and \cite[Corollary 1.9]{ARS} we have $\dim H_q(X; \F^{\QQ}_p) = h^{p,q}(\C V)$ and this completes the proof.
\end{proof}

Finally, the expression of the bounds in Theorem \ref{thm:main} is obtained by the  Lefschetz Hyperplane Theorem which implies that if  $p+q \neq n$ or $p \neq q$ then  then $h^{p, q}(\C V) = 0$. This completes the proof.

\begin{remark}
The above theorem also holds beyond the case of hypersurfaces. For example, Viro's patchworking construction has been generalised to complete intersections \cite{SturmfelsCompInt}, \cite{Bihan}. 
In order to obtain a version of Theorem \ref{thm:main} for complete intersections it remains to relate the dimensions of the tropical homology groups with $\Z_2$-coefficients to the Hodge numbers of complete intersections over the complex numbers. One possible route to making this connection is to establish tropical Lefschetz section theorems  and torsion freeness of the integral tropical homology groups for tropical  complete intersections as in the case for hypersurfaces in \cite{ARS}.

Even more generally, a tropical manifold is a polyhedral space locally modelled on matroid fans  \cite{MikRau}. A real phase structure on a tropical manifold would consist of  specifying  orientations of the local matroids, subject to compatibility conditions. It would be interesting to study real phase structures in this context and to generalise Theorem \ref{thm:toricvar}. 


\end{remark}

\subsection{Further consequences of the spectral sequence}

In Section \ref{sec:spectralsequence}, we describe how to go further in the spectral sequence to give criteria in  Theorem \ref{thmindivid}  for a real algebraic hypersurface arising from a primitive patchworking to attain the bounds in Theorem \ref{thm:main}. 
For a real algebraic variety, the Smith-Thom inequality  bounds the sum of the Betti numbers of the real part by the sum of the complexification, 
\begin{equation}\label{ST}
\sum_{q = 0}^n b_q(\R V)\leq \sum_{q = 0}^{2n} b_q(\C V).
\end{equation}
A real algebraic variety is called an \textit{$M$-variety}, or \textit{a maximal variety}, if  it satisfies equality in (\ref{ST}). 
The spectral sequence gives a necessary and sufficient condition for a  compact real algebraic hypersurface near the tropical limit to be maximal in the sense of the Smith-Thom inequality. 
Moreover,  Theorem \ref{thmindivid} gives a   criterion for  individual Betti numbers to  attain the bounds of Theorem \ref{thm:main}.

\begin{thm}\label{thmmax}

A compact real hypersurface with non-singular Newton polytope and near a non-singular tropical limit is maximal in the sense of the Smith-Thom inequality (\ref{ST}) if and only if the associated  spectral sequence  $(E^{\bullet}_{\bullet, \bullet}, \partial_{\bullet})$
degenerates at the first page.
\end{thm}
\begin{proof}
A compact real  hypersurface near the tropical limit is maximal in the sense of the Smith-Thom inequality if and only if all the inequalities in Theorem \ref{thm:main} are equalities. This happens exactly when the spectral sequence $(E^{\bullet}_{\bullet, \bullet}, \partial_{\bullet})$ degenerates at the first page. 
\end{proof} 

Viro proved the existence of non-singular maximal surfaces of any degree in $\PP^3$  \cite{Viro79}. Later, Itenberg and Viro proved that there exist non-singular projective hypersurfaces of any dimension that are asymptotically maximal \cite{IV}. This was generalised by Bertrand \cite{Bertrand06} to hypersurfaces and complete intersections in arbitrary toric varieties. Bertrand also proved in \cite{Bertrand06} that there exist toric varieties in any dimension that do not admit torically non-degenerate maximal hypersurfaces. However, all known examples are singular.
\begin{ques}
For every non-singular  Newton polytope $\Delta$, does there exist a maximal real hypersurface  in the toric variety  $Y_{\Delta}$ with Newton polytope $\Delta$?
\end{ques}

If the following conjecture is true, then it would simplify the statement of Theorem  \ref{thmindivid}, since we would only need to consider the differentials on the first page. 

\begin{conj}
For a compact non-singular real tropical hypersurface $(X,\mathcal{E})$  the spectral sequence associated to the filtration 
$$0 \subset C_{\bullet}(X ; \K_n)  \subset \dots \subset C_{\bullet}(X ; \K_1) \subset C_{\bullet}(X ; \mathcal{S}_{\mathcal{E}})$$
 degenerates at the second page. 
 \end{conj}

In Section \ref{sec:Haas}, we restrict our attention to real plane curves near the tropical limit. In this case, the only possible non-zero differential of the spectral sequence is on the first page. Using  the isomorphism  in Corollary \ref{cor:relativenadtrop}, this differential  is $\partial_1 \colon H_1(C; \F_0) \to H_0(C; \F_1)$. In Theorem \ref{thm:twisted}, we explicitly  describe this linear map using the twist description of patchworking for curves \cite[Section 3]{BIMS}. 
Using this description we also recover Haas' criterion for the maximality of curves  in toric surfaces near non-singular tropical limits in Theorem \ref{HaasMax}.

A real algebraic hypersurface near the non-singular  tropical limit has the same signature as the  Euler characteristic of its real part. This relation was first proved in the case of surfaces by Itenberg \cite{Itenberg1997},  and then later generalised to arbitrary  dimensions by Bertrand \cite{Bertrand10}. 
By  comparing  Euler characteristics of different pages of the above  spectral sequence, we recover this result and a generalisation to the non-compact case. 
A different proof of this signature formula  in the compact case was given by  Arnal in his master's thesis \cite{Arnal} also using tropical homology. 
We recall that as in Theorem \ref{thm:toricvar}, we let $Y_{\Delta}^o$ denote a partial compactification of the torus corresponding to a subfan of the dual fan of $\Delta$.

\begin{corollary}\label{cor:sign}
Let $V$ be a real algebraic hypersurface in the non-singular toric variety  $Y^o_{\Delta}$ with  Newton polytope $\Delta$ and near  the non-singular tropical limit $X$,   then 
$$\chi^{BM}(\R V)=\chi_{y}(\C V){{\big |}_{y = -1}},$$
 where $\chi_{y}(\C V)$ denotes the $\chi_y$ genus of $\C V$ and $\chi(\R V)$ 
is the topological Euler characteristic of $\R V$. 
In particular, if $Y_{\Delta^o}=Y_{\Delta}$ is compact, then
$$\chi(\R V)=\sigma(\C V),$$
where $\sigma(\C V)$ denotes the signature of $\C V$.

\end{corollary}

\begin{proof}

The Euler characteristic  of the $r$-th page of a spectral sequence is 
$$
\chi(E^r)=\sum_{p,q} (-1)^{q} \dim E^r_{q,p}.
$$
Since any page is by definition the homology of the preceding page, one has $$\chi(E^ \infty)=\cdots=\chi(E^0) =\sum_{p} \chi(C^{BM}_{\bullet}(X; \mathcal{F}_p)),$$
where $\chi(C^{BM}_{\bullet}(X; \mathcal{F}_p))$ denotes the Euler characteristic of the chain complex for the Borel-Moore homology of $X$ with coefficients in $\F_p$.

A complex hypersurface $\C V$ which is near a non-singular tropical limit is torically non-degenerate in the sense of \cite{DanilovKhovansky}. 
Therefore, by \cite[Theorem 1.8]{ARS}, the Euler characteristics of the chain complexes for tropical homology give the coefficients of the $\chi_y$ genus of $\C V$, 
$$\chi_y ( \C V) = \sum_{p}  (-1)^p \chi(C^{BM}_{\bullet}(X; \mathcal{F}_p))y^p. $$ Therefore, 
$$
\chi_y ( \C V){{\big |}_{y = -1}}=\chi(E^0)
$$
Moreover,  the Euler characteristic of the infinity page is equal to the  Borel-Moore Euler characteristic of $\R V$, 
$$
\chi(E^{\infty})=\chi^{BM}(\R V), 
$$
and this proves the first claim. 

Finally,  in the compact case since the  $\chi_y$ characteristic is defined with $(-1)^p$,  we  have $\sigma(\C V)=\chi_{-1}(\C V)$ and the corollary is proved.

\end{proof}

To end the introduction we would like to make a few remarks about the geometric inspiration behind the proofs of Theorem \ref{thm:main}. 
The construction of the  cosheaves $\S_{\mathcal{E}}$, $\K_p$, and $\F_p$, together with the linear maps $bv_p$, are all presented  using linear algebra, but their definitions are  geometrically  motivated. 
In the case of the sign cosheaf, the $\Z_2$-vector space $\S_{\mathcal{E}}(\sigma)$ associated to a face  $\sigma$ of $X$ is isomorphic to $H_{0}(\R P^{n} \backslash \R \A_{\sigma} ;\Z_2)$, where $\A_{\sigma}$ is a real hyperplane arrangement determined by the face $\sigma$ and $n$ is the dimension of $X$. 
The cosheaves $\F_p$ from tropical homology satisfy   $\F_p(\sigma) = H_{p}(\C P^{n} \backslash \C \A_{\sigma} ;\Z_2)$ \cite{Zharkov13}. This is described in the proof of Lemma  \ref{lem:dimFp}. 

For real varieties, the Viro homomorphism is a  partially defined multivalued homomorphism 
$$bv_{\ast} \colon H_{\ast}(\R V; \Z_2) \dashrightarrow  H_{\ast}(\C V; \Z_2),$$
where $H_{\ast}(\R V; \Z_2)$ and $H_{\ast}(\C V; \Z_2)$ denote the total homology of the real and complex parts respectively. 
A description of these homomorphisms is given in  \cite[ Appendix A2]{DegKhar}. 
The complement of a  real hyperplane arrangement $\A$ in $\R P^n$ is a disjoint  union of convex regions and therefore satisfies $H_{q}(\R P^n \backslash \R \A;  \Z_2) = 0$ for
all $q \neq 0$. Moreover, the complement of a real  hyperplane arrangement is a maximal variety in the sense of the Smith-Thom inequality (\ref{ST}) \cite[Introduction p.6]{OrlikTerao}. Therefore, in this special case the Viro homomorphism gives a collection of well defined graded maps 
$$bv_{p} \colon \Ker(bv_ {p-1})   \to   H_{p}(\C P^{n} \backslash \C \A ;\Z_2).$$
The map  $bv_0 \colon H_0(\R P^n \backslash \R \mathcal{A};\Z_2) \to H_{p}(\C P^{n} \backslash \C \A ;\Z_2)$ is induced by  the inclusion $i \colon \R P^n \backslash \R \A \to \C P^{n} \backslash \C \A $. To define the map $bv_p$, given  $\alpha\in \Ker(bv_{p-1})$, consider a $p$-chain $\beta$ in $\C P^n \backslash \C \A$ such that $\partial\beta=bv_{p-1}(\alpha)$. Then $bv_p(\alpha)$ is the homology class of the cycle $\beta+\mathrm{conj}(\beta)$. It follows from the maximality of $\C P^{n} \backslash \C \A$ that the complex conjugation acts as the identity on homology groups, see \cite[Corollary A.2]{Wilson}. Therefore, the maps $bv_p$ are well defined as they do not depend on the choice of the chain $\beta$.
Kalinin's spectral sequence \cite{Kalinin} induces a filtration on the real homology of a variety, which in the case of a real hyperplane  arrangement is given by  
$$ 0 = \Ker(bv_ {n})  \subset \Ker(bv_ {n-1})  \subset \dots \subset \dots \subset \Ker(bv_ {0})  \subset H_0(\R P^n \backslash \R \A; \Z_2).$$   

Although we do not use this geometry in the presentation of our arguments,  we borrow the notation for the Viro homomorphisms for our maps $bv_p \colon \K_p \to \F_p$ and use the letter $\K$ to denote the pieces of the 
filtration of $\mathcal{S}_{\mathcal{E}}$  in reference to Kalinin's filtration.

\subsection{Related works}
 
 In the case of real complete toric varieties, Hower \cite{Hower}  used a spectral sequence to relate the Betti numbers of a real toric variety to a $\Z_{2}$ variant of Brion's Hodge spaces for fans \cite{Brion}. The $\Z_2$-Hodge spaces for complete regular fans vanishes outside of the line $p=q$ \cite[Proof of Theorem 1.2]{Mtoric}, and nonsingular complete real toric varieties are all maximal in the sense of the Smith-Thom inequality. Hower proved that in the case of real toric varieties coming from reflexive polytopes, the spectral sequence also degenerates at the first page, and those toric varieties are again all maximal. Hower also exhibits an example of a six-dimensional projective toric variety which is not maximal, disproving a conjecture in \cite{Mtoric}. 

 The Hodge spaces for fans coincide with the tropical cohomology groups of the corresponding  tropical toric variety, since their defining chain complexes are isomorphic by definition. 
Finding necessary and sufficient conditions for  fans  to satisfy a version of Poincar\'e duality for Brion's Hodge spaces with integer coefficients would lead to a better understanding of fans defining maximal toric varieties.

 For Lagrangian toric fibrations equipped with an anti-symplectic involution,  Casta\~no-Bernard and Matessi study the cohomology of the fixed point locus using a long exact sequence which relates it to the cohomology of the Calabi-Yau manifold \cite{CastanoMatessi}. This is inspired by the Leray spectral sequence of Gross that relates the cohomologies of the Calabi-Yau manifold and the base space \cite{Gross}, see also  \cite{GrossSiebert}. Very recently Arg\"uz and Prince, computed the connecting maps of this long exact sequence  and the cohomology  groups with $\Z_2$-coefficients  of real Lagrangians in the quintic 3-fold and its mirror \cite{ArguzPrince}.    

\subsection{Outline of the paper}
 In Section \ref{section:tropicalhypersurface}, we review the definitions of tropical hypersurfaces  in toric varieties and their tropical homology groups. 
Section \ref{section:realtropicalhypersurfaces} uses  real phase structures on  non-singular tropical hypersurfaces  to describe Viro's primitive patchworking. In Subsection \ref{sec:signcosheaf},   we introduce the sign cosheaf and the  real tropical homology groups.  Section \ref{section:filtration} describes the augmentation filtration in general and adapts it to filter the sign cosheaves and the chain complexes. Section \ref{section:mainproof} contains the proof of Theorem \ref{thm:main}. Section \ref{sec:spectralsequence} investigates going further in the spectral sequence and lists all possible non-zero maps at further pages. 
Lastly in Section \ref{sec:Haas}, we illustrate the situation in the case of plane curves and describe the only possibly non-zero differential map in the spectral sequence in this case. This allows us to recover Haas' condition to obtain maximal curves.

\section*{Acknowledgement}
We are extremely grateful to Ilia Itenberg for his kind invitation to \'Ecole Normale Superieure and Paris VI and for the helpful discussions and suggestions.  
We would also like to thank Charles Arnal, Erwan Brugall\'e, Benoit Bertrand, Alfredo Hubard, Grisha Mikhalkin, Patrick Popescu-Pampu, Johannes Rau, Antoine Touz\'e, Jean-Yves Welschinger, and 
Ilia Zharkov for insightful discussions. 
We thank also H\"ulya Arg\"uz, Diego Matessi, Matilde Manzaroli and Thomas Prince for comments on a preliminary version of the paper.
This work was concluded while both authors were participants of the semester ``Tropical Geometry, Amoebas, and Polytopes" at the Institute  Mittag-Leffler. We are very grateful  to the organisers for their invitation and to the institute for their wonderful hospitality and working conditions. 

A.R. acknowledges support from the Labex CEMPI (ANR-11-LABX-0007-01). The research of  K.S. was supported by the Max Planck Institute Leipzig, and the Trond Mohn Foundation 
project ``Algebraic and topological cycles in tropical and complex geometry".

\section{Projective tropical hypersurfaces}
\label{section:tropicalhypersurface}
The tropical numbers are the set $\T = [-\infty, \infty)$. We equip $\T$ with the topology of a half open interval and $\T^n$ with the product topology. 
Tropical toric varieties are tropical manifolds in the sense of \cite{MikRau} with charts to $\T^{n}$. 
Just like toric varieties over a field, they are constructed 
 from rational polyhedral fans
see \cite[Section 6.2]{MacStu},  \cite[Section 3.2]{MikRau}. We recall that a  rational polyhedral fan $\Sigma$ is \emph{simplicial} if each of its cones is the cone over a simplex. A simplicial rational polyhedral fan is unimodular if the primitive integer directions of the rays of each cone can be completed to a  basis of $\Z^{n+1}$. 
A  tropical toric variety is \emph{non-singular} if it is built from a simplicial unimodular rational polyhedral fan. A tropical toric variety is compact if and only if the corresponding fan is complete.

A tropical toric variety $\T Y$ has a  stratification and the combinatorics of the stratification is governed by its fan $\Sigma$. 
A stratum of dimension $k$ of $\T Y$ corresponds to a cone $\rho$ of dimension $n+1-k$ of $\Sigma$. Denote the strata in  $\T Y$ corresponding to the cone $\rho$ simply by  $\T Y_{\rho}$.
 \begin{example}\label{ex:tropproj}
 The $(n+1)$ dimensional tropical torus is the space $\R^{n+1}$ and $(n+1)$-dimensional tropical affine space is $\T^{n+1}$.
  
 Tropical projective space  $\T P^{n+1}$  is constructed from the complete fan $\Sigma$ in $\R^{n+1}$ whose rays are in directions $-e_1, \dots -e_{n+1}, e_0 = \sum_{i = 1}^{n+1} e_i$, where $e_i$ denote the standard basis vectors. For every proper subset  $I \subset \{0, \dots n+1\}$ there is a cone of the fan defining $\T P^{n+1}$ of dimension equal to $|I|$. 
Moreover, analogous to projective space over a field, tropical projective space can also be defined as the quotient
$$\TP^{n+1}  =  \frac{\T^{n+2} \backslash (-\infty, \dots, -\infty)}{[x_0 : \dots: x_{n+1}] \sim [a+x_0: \dots: a+x_{n+1}]},$$
where $a \in \T \backslash -\infty$. 

Tropical projective space also admits a stratification determined by the fan $\Sigma$ defining it. 
The cones of $\Sigma$ correspond to a subsets $I \subsetneq  \{0, \dots, n+1\}$. We define the $I$-th open  stratum of $\TP^{n+1}$ to be
$$\TP_I^{n+1}   = \{ x \in \TP^{n+1} \ | \ x_i = -\infty \text{ iff } i \in I\}.$$
Notice that the open stratum $\TP_I^{n+1}$ can be identified with $\R^{n+1 - |I|}$. 
 \examp
\end{example}

When $\sigma$ is of dimension $k$ we have $\T Y_{\rho} \cong \R^{n+1-k}$, and $\T Y_{\rho}$ is a tropical torus of dimension $n+1-k$. 
For two cones $\rho$ and $\rho'$ of the fan  $\Sigma$ we have $\T Y_{\rho'} \subset \overline{\T Y_{\rho}}$ if and only if $\rho$ is a face of $\rho'$ in $\Sigma$. 
Morover if $\rho$ is a face of $\rho'$ in $\Sigma$, we have a projection map denoted by $\pi_{\rho,\rho'} \colon \T Y_{\rho}\rightarrow \T Y_{\rho'}$.
We assume the vertex of the fan to be $0$, so the corresponding open stratum of $\T Y$ is denoted by $\T Y_0$.

\subsection{Tropical hypersurfaces}

A tropical polynomial in $n+1$ variables is a function  $F_{\text{trop}} \colon \R^{n+1} \to \R$ of the form
\begin{equation}\label{tropPoly}
F_{\text{trop}}(x)=\max_{i \in A}(a_i+\langle i,x\rangle), 
\end{equation}
where $\langle \: \cdot \: , \cdot \: \rangle$ denotes the standard scalar product in $\R^{n+1}$, the set $A$ is a subset of  $\Z^{n+1}$, and $a_i\in\T$ for all $i\in A$.

A tropical polynomial of the form (\ref{tropPoly}) induces a regular subdivision of the  Newton polytope of its defining polynomial. 
A \emph{tropical hypersurface $X$ in $\R^{n+1}$} is the locus of non-linearity  of the function defined by a  tropical polynomial together  with weights naturally assigned to its top dimensional  faces, also known as \emph{facets}. 
The  tropical hypersurface of a polynomial  is dual to the regular subdivision of its Newton polytope induced by the convex-hull of the graph of $i\to a_i$, hence this subdivision is called the \emph{dual subdivision} of $X$.    
The weight of a facet is the integer length of the segment of the dual subdivision dual to the facet. We refer the reader to \cite[Section 3.1]{MacStu}, \cite[Section 5.1]{BIMS}, and  \cite[Section 2.3]{MikRau} for further details and  examples. 

A tropical hypersurface in $\R^{n+1}$ is \emph{non-singular} if its dual subdivision is primitive, meaning that each $n+1$ dimensional polytope of the subdivision has normalised lattice volume equal to $1$. In particular, the  weights on all facets of a non-singular tropical hypersurface are  equal to one.
We define a tropical hypersurface in a tropical toric variety $\T Y$ to  be the closure of a tropical hypersurface in $\R^{n+1} = \T Y_0 \subset \T Y $ where $\T Y_0 \cong \R^{n+1}$ is the open stratum corresponding to the vertex $0$ of the fan defining $\T Y$. 

Given a tropical hypersurface $X$ in $\R^{n+1}$ with Newton polytope $\Delta$, we will consider its (partial) compactification in a tropical toric variety $\T Y_{\Delta}^o$ which is defined by a subfan of the dual fan of $\Delta$. In this case, the compactification of $X$ in $\T Y_{\Delta}^o$  is a non-singular tropical variety if  $\T Y_{\Delta}^o$ is non-singular, which is guaranteed if Newton polytope $\Delta$ is non-singular.

\begin{example}\label{ex:hypersurfaceinTP}
If $X$ is a non-singular  tropical hypersurface in $\T P^{n+1}$, then the Newton polytope of $X$ is equal to $d\Delta_{n+1}$ for some $d$, where 
$$d\Delta_{n+1} = \text{ConvHull}\{0, de_1, \dots, de_{n+1}\}. $$
If $X$ is non-singular  in $\R^{n+1}$ with the above Newton polytope, then it is dual to a unimodular subdivision of $d\Delta_{n+1} $. The intersection  $X_I \subset \TP^{n+1}_I$ is also a  non-singular tropical hypersurface dual to the subdivision of the corresponding face of $d \Delta_{n+1}$.
\end{example}

Our convention is that all  faces of a hypersurface $X$ in a tropical toric variety $\T Y$ are  closed. We let $X_{\rho} := X \cap \T Y_{\rho}$. The faces of $X_{\rho}$  are also considered  closed in $\T Y_{\rho} \simeq \R^{n+1- \dim \rho}$. 
For a face $\sigma$ of a tropical hypersurface,  we let $\relint(\sigma)$ denote its relative interior. 
The \emph{sedentarity} of a point $y$ in $\T Y$ is $\rho$ if $y$ is contained in the stratum $\T Y_{\rho}$. 
The \emph{sedentarity} of a face $\sigma$ is  denoted by   $\sed(\sigma)$ and is equal to $\rho$ if $\relint(\sigma) \subset \T Y_\rho$. 
The \emph{parent face} of a face $\tau$ of $X$ of dimension $k$ and sedentarity $\rho$ is the unique face $\sigma$ of $X$ of empty sedentarity and dimension $k+\dim \rho$ such that $\tau$ is in the boundary of $\sigma$. The \emph{star} of a face $\sigma$ is is the star of any of its relative interior point 
$$
\mathrm{Star}(\sigma)=\left\lbrace v\in\R^{n+1-|\sed(\sigma)|} \mid \exists \epsilon >0, x_\sigma+\epsilon v \in \sigma \right\rbrace,
$$
where $x_\sigma$ is any point in $\relint(\sigma)$.

\subsection{Tropical homology}
\label{section:tropicalhomology}

The cosheaves that we use throughout the text will always be vector spaces over $\Z_2$. 
Let $X$ be a non-singular tropical hypersurface in a tropical toric variety $\T Y$. Let the defining fan of  $\T Y$ be the simplicial unimodular fan $\Sigma$ in $\R^{n+1}$, and let $\Z^{n+1}$ denote the standard lattice in $\R^{n+1}$. 
Let  $\rho$ be a cone of $\Sigma$ of dimension $s$ with primitive integer generators  $r_1, \dots, r_s$ and define the integral tangent space to $\T Y_{\rho}$ as
$$T_{\Z}(\T Y_{\rho}) := \frac{\Z^{n+1}}{ \langle r_1, \dots , r_s \rangle}.$$
For a face $\sigma$ of $X$ of sedentarity $\rho \in \Sigma$,  
 let $T_ \Z(\sigma) \subset T_ \Z(\T Y_{\rho})$ denote the integral tangent space of $\relint(\sigma)$. 
  Since $X$ is a non-singular hypersurface the reduction modulo $2$ of the free  $\Z$-module $T_{\Z}(\sigma)$ is a vector space of the same dimension. In fact, at any vertex $v$ adjacent to $\sigma$, one can complete a basis of the free $\Z$-module $T_{\Z}(\sigma)$ into a basis of $\Z^{n+1-\dim \rho}$ with vectors in $T_{\Z}(\sigma_i)$ for $\sigma_i$ faces adjacent to $v$. We denote this vector space over $\Z_2$ by  $\F_1(\sigma)$. 

If  $\T Y_{\rho}$ and $\T Y_{\eta}$ are a pair of strata corresponding to cones $\rho \subset \eta$ of $\Sigma$  then  $\T Y_{\eta} \subset \overline{\T Y}_{\rho}$ then the generators of the cone $\eta$ contain the generators of the cone $\rho$. Therefore there is a   projection map: 
\begin{equation}\label{eq:projmaps}
\pi_{\rho \eta} \colon T_ \Z (\T Y_{\rho})  \to T_ \Z (\T Y_{\eta}). 
\end{equation}
Upon taking the reduction modulo $2$ we get a map $\pi_{\sigma \tau} \colon \F_1 (\T Y_{\rho})  \to \F_1 (\T Y_{\eta})$. For faces $\sigma$ and $\tau$ in $\T Y_{\rho}$  and $\T Y_ \eta$, respectively produces a map which we denote 
\begin{equation}\label{projmap}
\pi_{\sigma \tau} \colon \F_1 (\sigma)  \to \F_1 (\tau).  
\end{equation}

\begin{definition}
Let $X$ be a non-singular  tropical hypersurface in a tropical toric variety $\T Y$ with defining fan $\Sigma$. 
The $p$-multi-tangent spaces of $X$ are cellular cosheaves $\F_p$ on $X$. For $\rho \in \Sigma$ and $\tau$ a face  of $X_{\rho} := X \cap \T Y_{\rho}$  we have 
\begin{equation}\label{def:Fp}
\F_p(\tau) = \sum_{\substack{\tau \subset \sigma \subset X_I \\ \dim(\sigma) = n- \dim \rho}} \bigwedge^p \F_1(\sigma). 
\end{equation}

When  $\tau\subset\sigma$, the maps of the cellular cosheaf 
$i_{\sigma \tau} \colon \F_p(\sigma) \to \F_p(\tau)$ 
are induced by the inclusions  $\F_1(\sigma) \to \F_1(\tau)$ when $\sigma$ and $\tau$ have the same sedentarity and otherwise are induced by the quotient map $\pi_{\sigma \tau}$
from (\ref{projmap}). 
\end{definition}

\begin{example}\label{ex:FpPlane}
The tropical plane $X \subset \T P^3$ is the closure of a two dimensional  fan $X_{0}$ in $\R^3$. The fan $X_{0}$ has rays $\tau_1,  \tau_2, \tau_3, $ and $\tau_0$  in respective directions $-e_1, -e_2, -e_3,$ and $e_0 = e_1+ e_2+e_3$, where $e_i$'s are the standard basis vectors. Every  pair of rays generates a two dimensional face of $X \cap \R^3$,    see the right hand side of Figure \ref{figplane}. 
Denote by $\va_i$ the reduction of $e_i$ mod $2$. 

Let $\sigma_{ij}$ denote the two dimensional face spanned by rays $\tau_i$ and $\tau_j$. 
Then $\F_1(\sigma_{ij}) = \Z_2 \langle \varepsilon_i, \varepsilon_j \rangle$ and $\F_2(\sigma_{ij}) = \Z_2 \langle \varepsilon_i \wedge \varepsilon_j \rangle$.  For the ray $\tau_i$, we
obtain $\F_1(\tau_i) = \Z_2^3$ and $\F_1(\tau_i) = \Z_2 \langle \varepsilon_i \wedge \varepsilon_j, \varepsilon_i \wedge \varepsilon_{j'} \rangle$ where $\varepsilon_i, \varepsilon_j, \varepsilon_{j'}$ form a basis of $\Z_2^3$. 
For any face $\tau$ of $X$, we have  $\F_0(\tau) = \Z_2$.  
\end{example}

\begin{lemma}\label{lem:dimFp}

Let $X$ be an $n$-dimensional  non-singular tropical hypersurface of a tropical toric variety $\T Y$. 
For a face $\tau$ of $X$ of dimension $k$ and  sedentarity $\rho$  the polynomial defined as 
$$\chi_{\tau}(\lambda) := \sum_{p = 0}^{n} (-1)^p\dim \F_p(\tau) \lambda^p,$$
is equal to 

$$
\chi_{\tau}(\lambda)=(1-\lambda)^{k}\left[(1-\lambda)^{n-k+1-\dim \rho}-(-\lambda)^{n-k+1-\dim \rho}\right]. 
$$

\end{lemma}

\begin{proof}
By \cite[Theorem 4]{Zharkov13}, the $\mathbb{Z}$-multi-tangent spaces  $\F_p^{\mathbb{Z}}(\sigma)$ are isomorphic to the dual of the  $p$-th graded piece of the Orlik-Solomon algebra of the matroid of  an associated projective  hyperplane   arrangement $\A_{\sigma}$ defined over the complex numbers. The Orlik-Solomon algebra of this arrangement is isomorphic to the cohomology ring of the complement of the arrangement in projective space so that $\text{Hom}(\F_p^{\Z}(\sigma), \Z) \cong H^p(\C P^n \backslash \C \A_{\sigma};\Z).$ The homology groups of the complement of a complex hyperplane arrangement are torsion free so $\F_p(\sigma) \cong H_p(\C P^n \backslash \C \A_{\sigma};\Z_2).$

Let $\mathcal{P}_{n-k}$ denote the  $(n-k)$-dimensional pair of pants;  that is the complement of $n-k+2$ hyperplanes in general position in $\C P^{n-k}$. 
For a  face $\tau$ of $X$ of dimension $k$ and sedentarity $\rho$, the complement of the associated arrangement is $\mathcal{P}_{n-k-\dim \rho} \times (\C^*)^k$. 
Therefore,  we have the isomorphism $\F_p(\tau) \cong H_p(\mathcal{P}_{n-k-\dim \rho} \times (\C^*)^k ; \Z_2)$,   and by  the K\"unneth formula for the homology groups we have 
$$
\chi_{\tau}(\lambda) = \chi_{\mathcal{P}_{n-k-\dim \rho} }(\lambda) \chi_{(\C^*)^k}(\lambda), 
$$
where $\chi_{\mathcal{P}_{n-k-\dim \rho} }(\lambda) $ and $ \chi_{(\C^*)^k}(\lambda)$ are the Euler-Poincar\'e polynomials of $\mathcal{P}_{n-k-\dim \rho} $ and $(\C^*)^k$ respectively. 
Calculating the homology of these spaces shows that 
\begin{align*}
\chi_{\mathcal{P}_{n-k} }(\lambda) & = \sum_{r=0}^{n-k-\dim \rho}(-1)^r\binom{n-k-\dim \rho+1}{r}\lambda^r \\
&= (1-\lambda)^{n-k+1-\dim \rho}-(-\lambda)^{n-k+1-\dim \rho},
\end{align*}
and 
$$\chi_{(\C^*)^k}(\lambda) = \sum_{s = 0}^k (-1)^s \binom{k}{s} \lambda^s = (1-\lambda)^{k}.$$
The product of these two polynomials is precisely the description of $\chi_{\tau}(\lambda)$ in the lemma. 
\end{proof}

\begin{definition}\label{def:trophomo}
Let $X$ be a non-sinuglar tropical hypersurface of a tropical toric variety $\T Y$. 
The groups of \emph{cellular $q$-chains with coefficients in $\mathcal{F}_p$} are 
$$C_q(X ; \mathcal{F}_p) = \bigoplus_{\dim \sigma = q} \mathcal{F}_p(\sigma).$$
The boundary maps $\partial \colon C_q(X ; \mathcal{F}_p) \to C_{q-1}(X ; \mathcal{F}_p) $ are given by the direct sums of the cosheaf maps $i_{\sigma \tau}$ for $\tau\subset \sigma$. 
The \emph{$(p,q)$-th tropical homology group} is 
$$H_q(X; \mathcal{F}_p) = H_q(C_{\bullet}(X ;\mathcal{F}_p)).$$

\end{definition}

\section{Real tropical hypersurfaces}
\label{section:realtropicalhypersurfaces}

In this section we describe a tropical approach to Viro's primitive patchworking construction via real phase structures. For an explanation of how it relates to Viro's original construction see Remark \ref{Viro}. 
Section \ref{realphasesub} defines real phase structures on tropical  hypersurfaces and describes how to obtain the real part of the tropical hypersurface. 
In Subsection \ref{sec:signcosheaf}, we introduce the sign cosheaf on a tropical hypersurface and prove that its homology groups are isomorphic to the homology groups of the real part of a tropical variety equipped with a real phase structure.

\subsection{Real phase structures and patchworking}
\label{realphasesub}
\begin{defi}

\label{realphase}
A \emph{real phase structure} on an $n$-dimensional  non-singular tropical hypersurface $X$ in a tropical toric variety $\T Y$ is a collection $\mathcal{E} = \{\mathcal{E}_{\sigma} \}_{\sigma \in \text{Facet}(X_0)^n}$  where   $\mathcal{E}_\sigma \subset \Z_2^{n+1}$ is an $n$-dimensional  affine subspace
parallel 

to $\F_1(\sigma)$. 
The collection $\mathcal{E}$ must satisfy the following property:
 
{\center{If $\tau$ is a face of $X \cap \T Y_0$ is  of codimension $1$, then for any facet $\sigma$ adjacent to $\tau$ and any element $\varepsilon\in\mathcal{E}_\sigma$, there exists a unique facet $\sigma'\neq \sigma$ adjacent to $\tau$ such that $\varepsilon\in\mathcal{E}_{\sigma'}$.}}

A non-singular tropical hypersurface equipped with a real phase structure is called a \emph{non-singular real tropical hypersurface}.

\end{defi}

\begin{exa}
\label{exrealtropline}
 Figure \ref{realtropline} depicts a real tropical line $X$ in the tropical projective plane $\TP^2$. On each edge $\sigma_0, \sigma_1, \sigma_2$ of the line there is a set  of vectors in $\Z_2^{n+1}$. These vectors indicate all the points in the affine subspace $\mathcal{E}_{\sigma_i}$ for a real phase structure $\mathcal{E}$. 
 
 The vertex of the tropical line is the only codimension one face.  For $(0,0) \in \mathcal{E}_{\sigma_1}$,  we have that  $(0,0) \in \mathcal{E}_{\sigma_2}$ and  $(0,0) \not \in \mathcal{E}_{\sigma_0}$. This is the condition in Definition \ref{realphase} for the face $\sigma_1$ and the element $(0,0)$. 
 
\begin{figure}
 \begin{minipage}[l]{.46\linewidth}
  \centering
  \includegraphics[width=5cm,height=4.5cm]{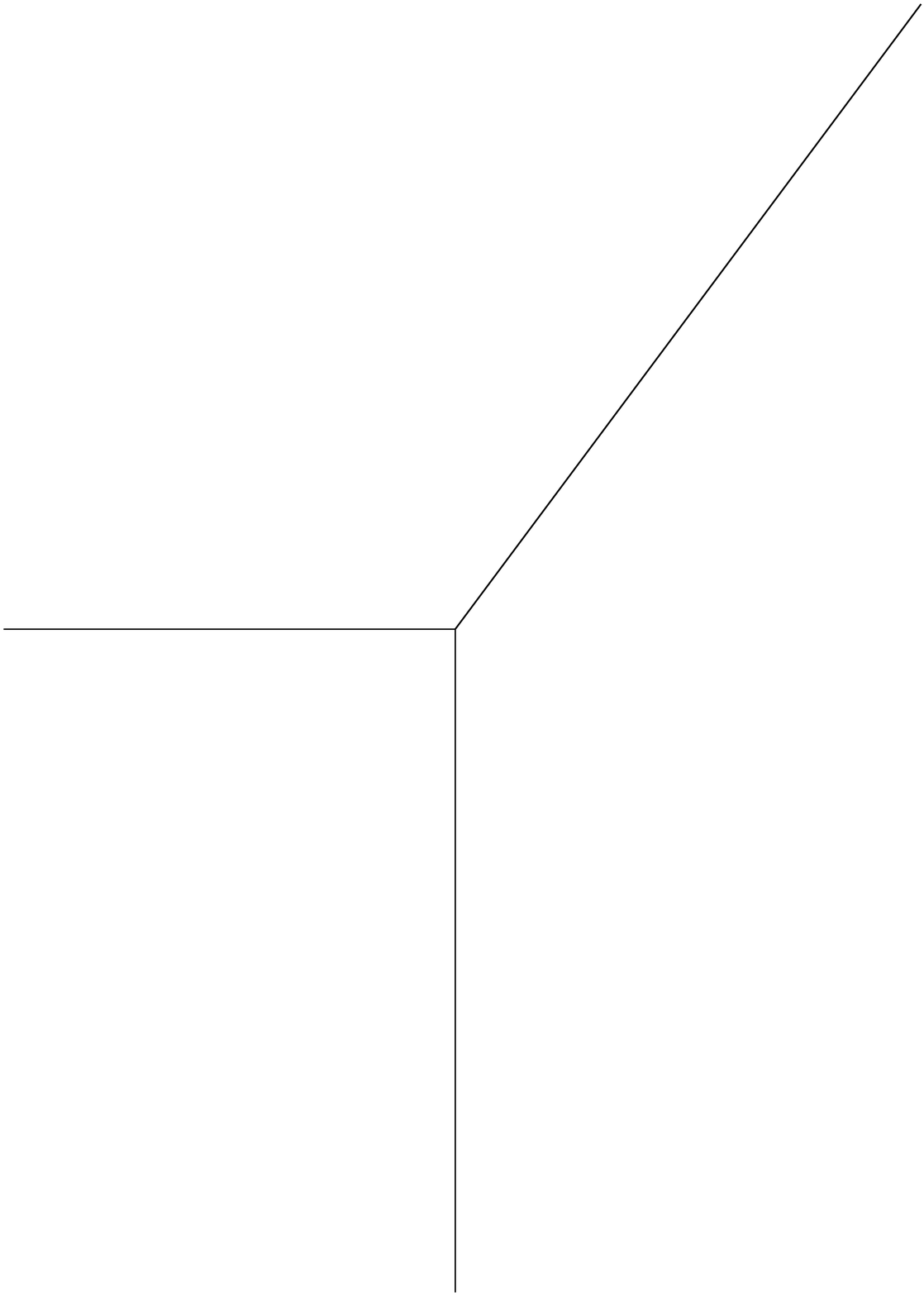}
  \put(-140,70){\tiny $\left\lbrace (0,0), (1,0) \right\rbrace$} 
  \put(-70,20){\tiny $\left\lbrace (0,0), (0,1) \right\rbrace$}
  \put(-15,110){\tiny $\left\lbrace (0,1), (1,0) \right\rbrace$}
 \end{minipage} \hfill
 \begin{minipage}[l]{.46\linewidth}
  \centering \includegraphics[width=5cm,height=5cm]{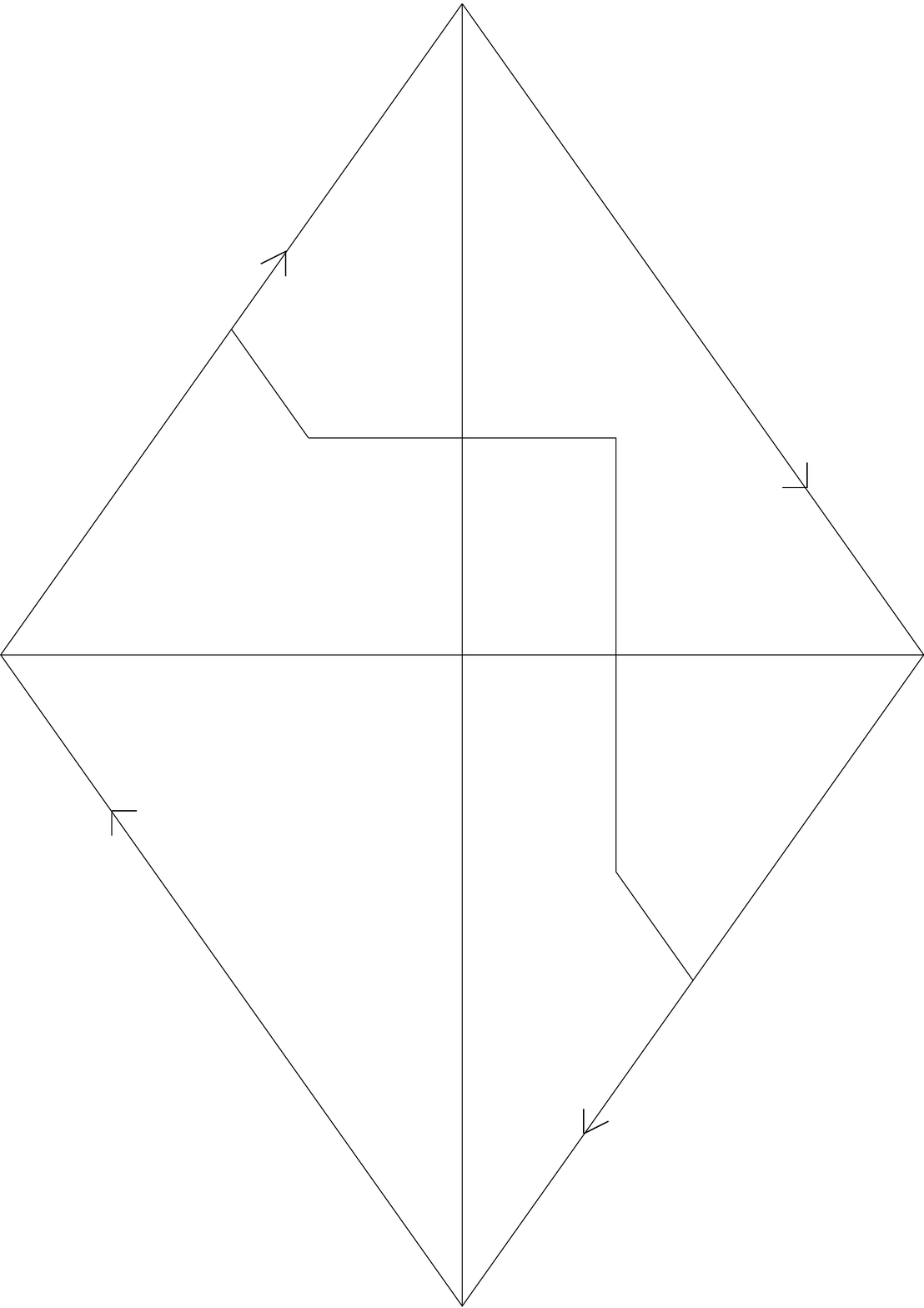}
 \end{minipage}
 \caption{On the left is the real tropical line $X \subset \TP^{2}$ with a real phase structure $\mathcal{E}$ from Example \ref{exrealtropline}. On the right hand side is its real part $\R X_{\mathcal{E}}$ in $\R P^{2}$. \label{realtropline}} 
\end{figure}
\end{exa}

\begin{exa}
\label{exrealtropplane}
Recall the tropical hyperplane described in Example \ref{ex:FpPlane} and depicted on the right of  Figure \ref{figplane}.  

The following collection of  affine spaces forms a real phase structure on $X$, %
$$\mathcal{E}_{\sigma_{12}}  = \langle \va_1, \va_2 \rangle, \quad \mathcal{E}_{\sigma_{13}}  = \langle \va_1, \va_3 \rangle, \quad \mathcal{E}_{\sigma_{23}}  = \langle \va_2, \va_3 \rangle,$$
$$\mathcal{E}_{\sigma_{01}}  = \langle \va_0, \va_1 \rangle + \va_3, \quad \mathcal{E}_{\sigma_{02}}  = \langle \va_0, \va_2 \rangle + \va_1, \text{   and    }  \mathcal{E}_{\sigma_{03}}  = \langle \va_0, \va_3 \rangle + \va_2. $$

Given a plane $P \subset \mathbb{P}^3$ defined over the real numbers, the intersection of $L$ with the coordinate hyperplanes of $\mathbb{P}^3$ defines an arrangement of real hyperplanes on $\R P \cong \R P^2$. Such is the picture on the left hand side of Figure \ref{figplane}.  Each region of the complement of this hyperplane arrangement on $\R P \subset \R P^3$ lives in a single orthant of $\R^3  = \R P^3 \backslash \{x_0 = 0\}$. In Figure \ref{figplane}, each connected component of the complement of the line arrangement is labelled with the vector in $\Z_2^{3}$ corresponding to this orthant.  
Let $L_i = \{x_i=0\} \cap P \subset \R P^3$ and set  $p_{ij} = L_i \cap L_j$. 
Notice that the points contained in the affine space $\mathcal{E}_{\sigma_{ij}}$ of the real phase structure on $X$ coincide with the collection of signs of the regions of the complement of the line arrangement which are adjacent to the point $p_{ij}$. 

\begin{figure}
 \begin{minipage}[l]{.46\linewidth}
  \centering
\begin{tikzpicture} 
\draw (-2.25,0) -- (0.25,4) 
node[anchor=south west] {$L_2$} 
node[pos=0.6, right]{\tiny{$(1, 0, 0)$}} 
node[pos=0.325, right]{\tiny{$(0, 0, 0)$}}
node[pos=0.25, left]{\tiny{$(0, 1, 0)$}} 
node[pos=0.75, left]{\tiny{$(1, 1, 0)$}};
\draw (2.25,0) -- (-0.25,4)node[anchor=south east] {$L_{\infty}$} 
node[pos=0.325, left]{\tiny{$(1, 0, 1)$}} node[pos=0.25, right]{\tiny{$(0, 1, 0)$}} node[pos=0.75, right]{\tiny{$(0, 1, 1)$}};
\draw (-2.75,0.) -- (2.5,2.5)node[anchor=south west] {$L_3$}
node[pos=0.37, right]{\tiny{$(0, 0, 1)$}} ;
\draw (2.75,0.) -- (-2.5,2.5)node[anchor=south east] {$L_{1}$};
  \end{tikzpicture}
 \end{minipage} \hfill
 \begin{minipage}[l]{.46\linewidth}
  \centering \includegraphics[width=6cm,height=5cm]{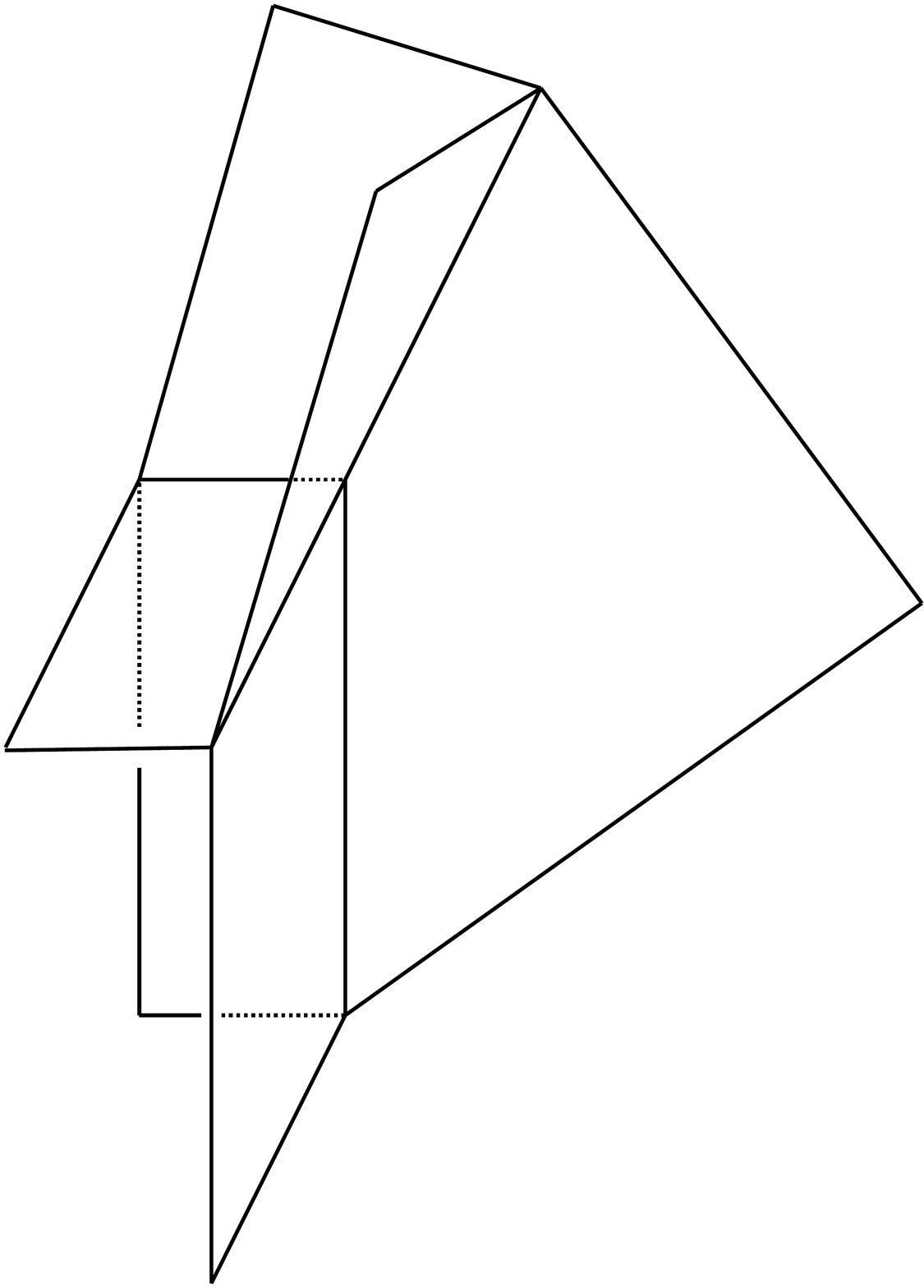}
 \put(-160,65){$\sigma_{12}$}
  \put(-102,113){$\sigma_{01}$}
   \put(-50,65){$\sigma_{03}$}
    \put(-120,130){$\sigma_{02}$}
     \put(-130,25){$\sigma_{13}$}
      \put(-160, 35){$\sigma_{23}$}
 \end{minipage}

\caption{The left hand side depicts a real line arrangement in $\R P^2$ arising from a linear embedding $\R P^2 \to \R P^3$. On the right is a  tropical plane $X$ in $\T P^3$. A real phase structure  on $X$ is described in Example \ref{exrealtropplane}.}\label{figplane}
\end{figure}
\end{exa}

Following \cite[Chapter 11]{GKZ}, we now describe how to obtain a space homeomorphic to 
the real part of a toric variety by glueing together multiple copies of a tropical toric variety. 
Let $\Sigma$ be a rational  polyhedral fan in $\R^{n+1}$ defining a tropical toric variety $\T Y$. 
For every $\varepsilon \in \Z_2^{n+1}$, let $\T Y(\varepsilon)$ denote a copy of $\T Y$ indexed by $\varepsilon$. 
Define 
\begin{equation}\label{eq:tropicalrealtoric}
\tilde{\T Y} :=  \bigsqcup_{\varepsilon \in \Z_2^{n+1}} \T Y({\varepsilon})/ \sim,  
\end{equation}
where $\sim$ identifies strata $\T Y_{\rho}(\varepsilon)$ and $\T Y_{\rho}(\varepsilon')$ if and only if $\varepsilon +  \varepsilon'$ is in the reduction modulo $2$ of the linear space spanned by the cone $\rho$. 

The following theorem is a direct translation of  \cite[Theorem 11.5.4]{GKZ}.
\begin{thm} 
The topological space $\tilde{\T Y}$ is homeomorphic to the real point set of the toric variety $\R Y$.

\end{thm}

\begin{example}\label{ex:realprojspace}
We we explicitly describe how to obtain $\R P^{n+1}$ by glueing together multiple copies of $\TP^{n+1}$.
For every $\varepsilon=(\va_1,\cdots ,\va_{n+1})\in \Z_2^{n+1}$, let $\TP^{n+1}{(\varepsilon)}$ denote a copy of 
$\TP^{n+1}$ indexed by $\varepsilon$. Then

$$\R P^{n+1} \cong \bigsqcup_{\varepsilon \in \Z_2^{n+1}} \TP^{n+1}({\varepsilon})/ \sim,  $$
where $\sim$ is the equivalence relation generated by identifying $ x \in \TP^{n+1}({\varepsilon}) $ and $ x' \in \TP^{n+1}({\varepsilon'}) $ for $\varepsilon \neq \varepsilon'$, such that 
$\left[x_0,\cdots,x_{n+1}\right]= \left[x'_0,\cdots,x'_{n+1}\right]$ 

and

\begin{itemize}
\item if $x_0 \neq -\infty$, then  there exist a unique $1\leq j \leq n+1$ such that $\varepsilon_j \neq \varepsilon'_j$. Moreover, we must have $x_j = x_j' = -\infty$. 
\item if $x_0 = -\infty$,  then we must have $\va_i  \neq \va_i^{\prime}$ for all $1 \leq i \leq n+1$.

\end{itemize}

\end{example}

Given a polyhedron $\sigma$  of sedentarity $0$ contained in $\T Y$ and $\varepsilon \in \Z_2^{n+1}$ we let $\sigma^{\varepsilon}$ denote its   copy in $\T Y{(\varepsilon)}$.

\begin{defi}
Let $X$ be a non-singular tropical hypersurface in  a tropical toric variety $\T Y$ together with a real phase structure $\mathcal{E}$. 

The \emph{real part of $X$} with respect to the real phase structure $\mathcal{E}$ is denoted $\R X_{\mathcal{E}}$ and is the  image in $\R Y$ of 
$$
\bigcup_{\substack{\mbox{ \tiny{facets of} } \sigma\subset X_0 \\ \varepsilon\in\mathcal{E}_\sigma}} \overline{\sigma}^{\varepsilon}.
$$
where $\overline{\sigma}^{\varepsilon}$ denotes the closure of $\sigma^{\varepsilon} \subset \T Y_0(\varepsilon)$  in $\T Y ({\varepsilon})$. 

\end{defi}

The following theorem is the tropical reformulation of a particular case of the combinatorial version of Viro's patchworking theorem from  \cite{Viro84}.

\begin{thm}[Viro's patchworking \cite{Viro84}]\label{Patch}
Let $(X,\mathcal{E})$ be a non-singular real tropical hypersurface with Newton polytope $\Delta$ in a non-singular tropical toric variety $\T Y_{\Delta}^o$ corresponding to a subfan of the dual fan of $\Delta$. 
Then there exists a non-singular real algebraic hypersurface $V$ of  $Y_{\Delta}^o$ also with Newton polytope $\Delta$ such that $$
(\tilde{\T Y}_{\Delta}^o, \R X_\mathcal{E} ) \cong (\R Y_{\Delta}^o ,\R V).
$$
\end{thm}

\begin{remark}
\label{Viro}
For the reader's convenience, we explain the connection between the tropical version of primitive patchworking and Viro's original formulation  as described in \cite{Viro84, Itenberg93,JJR}. 

The input of Viro's original formulation of primitive  patchworking is a  regular  subdivision of a lattice polytope $\Delta \subset \R^{n+1}$, whose normal fan is unimodular,  
 together with a choice of sign  $\delta_{i} \in \{ + , - \}$ for each lattice point $i \in \Delta \cap \Z^{n+1}$. 

Given a tropical hypersurface $X \in \R^{n+1}$, its dual subdivision is a regular subdivision of $\Delta$ which, by definition, is primitive if $X$ is non-singular. 
Every edge $e$ of the dual subdivision of $\Delta$ is dual to a facet $\sigma_e$ of $X$. 
From a  real phase structure $\mathcal{E}$ on $X$ we produce a  collection of signs $\delta_i$ for all $i \in \Delta \cap \Z^{n+1}$ as follows. 
Two vertices of an edge $e$ of the subdivision of $\Delta$ are assigned different signs if and only if  $\mathcal{E}_{\sigma_e}$ contains the origin $(0,\cdots,0)$ in $\Z_2^{n+1}$. For more details we refer to \cite[Lemma 1]{Renau}.
Upon choosing the sign of one lattice point in $\Delta$ arbitrarily, this rule determines a collection of signs for each integer point in $\Delta$.  

From the subdivision of $\Delta$ and the assignment of signs to all lattice points in $\Delta$, Viro's construction builds a polyhedral complex in the following way. 
For $\varepsilon \in \Z_2^{n+1}$, let $\Delta(\varepsilon)$ denote the symmetric copy of $\Delta$ in the orthant of $\R^{n+1}$ corresponding to $\varepsilon$. 
Then define 
\begin{equation}\label{RPnagain}
\tilde{\Delta} :=  \bigcup_{\varepsilon\in\Z_2^{n+1}}\Delta(\varepsilon) / \sim, 
\end{equation}
where the equivalence relation $\sim$ is the same as described for the tropical toric variety coming from the fan $\Sigma$ which is a subfan of the dual fan of $\Delta$.

The  triangulation of $\Delta$ induced by $X$ induces a symmetric triangulation of $\Delta(\varepsilon)$. Moreover, the sign choices $\delta_i \in \{+, -\}$ for $i \in \Delta \cap \Z^{n+1}$ induce choices of signs for $\Delta(\varepsilon)$ for all $\varepsilon \in \Z_2^{n+1}$ by way of the following rule: 
For $i_1, \dots, i_n \in \Delta(\varepsilon)$ 
$$\delta_{(i_1, \dots, i_n) }= \left(\prod_{j=1}^{n+1}(-1)^{\varepsilon_j i_j}\right)\delta_{(| i_1|,...,| i_{n+1} |  )}.$$

In other words, when passing from a lattice point to its reflection in a coordinate hyperplane, the sign is preserved if the distance from the lattice point to the hyperplane is even, and the sign is changed  if the distance is odd.

For a  simplex $T$ in the subdivision  of $\Delta(\varepsilon)$ let  $S_T$ denote the convex hull of the midpoints of the edges of $T$ having endpoints of opposite signs. Denote by $S$ the union of all such $S_T$ considered in the quotient to $\tilde{\Delta}$ as in (\ref{RPnagain}). Then $S$ is an $n$-dimensional  piecewise-linear manifold contained in $\tilde{\Delta}$. It turns out that pairs $(\tilde{\T Y} ,\R X_{\mathcal{E}})$ and $(\tilde{\Delta} ,{S})$ are combinatorially isomorphic and homeomorphic. Thus the two formulations of patchworking are equivalent. For more details see \cite[Lemma 1]{Renau}.

From here a  polynomial defining the  hypersurface $V$ from Theorem \ref{Patch} can be written down explicitly. The defining polynomial of $V$ is 
\begin{equation}
{\bf F}_t(x)=\sum_{(i_1,...,i_{n+1})\in \Z^{n+1}\cap \Delta}\delta_{i_1,...,i_n}\left(x_1^{i_1}\cdots x_{n+1}^{i_{n+1}} x_0^{d-\sum i_j}\right)
t^{-a_{(i_1,...,i_n)}},
\label{Viropoly}
\end{equation}
where the $a_{(i_1,...,i_n)}$'s are the coefficients from the tropical polynomial in $(\ref{tropPoly})$ and $t>0$ is a  sufficiently large real number.

\end{remark}

\begin{definition}\label{def:neartroplim}
A real algebraic hypersurface $V$ in a  toric variety $Y$ 
 is called \emph{near a non-singular tropical limit} if it is defined by a polynomial  ${\bf F}_t(x)$   of the form $(\ref{Viropoly})$ for $t$ sufficiently large coming from a non-singular tropical hypersurface $X$ with a real phase structure $\mathcal{E}$ and the fan defining $Y$ is a subfan of the fan dual to the Newton polytope of ${\bf F}_t$ 

\end{definition}

In particular, a  hypersurface near a non-singular tropical limit $X$ with real phase structure ${\mathcal{E}}$  will satisfy the homeomorphism of pairs from Theorem \ref{Patch} 
$$(\tilde{\T Y}_{\Delta}^o,\R X_\mathcal{E})\cong (\R Y_{\Delta}^o ,\R V).$$

\subsection{The sign cosheaf}
\label{sec:signcosheaf}

Let $X$  be a non-singular real tropical hypersurface equipped with a real phase structure $\mathcal{E}$. 
By definition, for any facet $\sigma$ of $X$ of sedentarity $0$,  the real phase structure  $\mathcal{E}$ gives an affine space $\mathcal{E}_\sigma$ of direction $\F_1(\sigma)$. Let us extend the real phase structure to facets of higher sedentarity as follows. 

Recall the definition of the map $\pi_{\sigma \tau} \colon \F_1(\sigma) \to \F_1(\tau)$ from (\ref{projmap})  via the   projection maps $\pi_{\rho \eta} \colon \F_1(\T Y_{\rho} )\to \F_1(\T Y_\eta)$ when $\sigma$ and $\tau $ are in strata $\T Y_{\rho}$ and $Y_{\eta}$, respectively. 
Also recall that the parent face of a face $\tau$ of $X$ of dimension $k$ and sedentarity $\rho$ is the unique face $\sigma$ of $X$ of empty sedentarity and dimension $k+\dim \rho$ such that $\tau$ is in the boundary of $\sigma$.
Let  $\tau_0$ denote the parent face of a facet $\tau$  in $X_\rho$ so that  $\pi_{\tau_0 \tau}(\tau_{0})=\tau$. Define $\mathcal{E}_{\tau}=\pi_{\tau_{0} \tau}(\mathcal{E}_{\tau_{0}})$. Notice that $\mathcal{E}_{\tau}$ is an affine space of $\F_1(\T Y_ \rho)$ which is parallel to $\F_1(\tau)$.

\begin{exa} \label{exreallineext}
The tropical line in $\T P^2$ from Example \ref{exrealtropline} contains three points of non-empty sedentarity. The projection of the affine vector spaces $\mathcal{E}_{\sigma_1},
\mathcal{E}_{\sigma_2}$ for the horizontal and vertical edges are simply $0 \in  \F_1(\TP^{n+1}_{\{1\}})$ and  $0 \in  \F_1(\TP^{n+1}_{\{2\}})$, respectively. For the  diagonal edge  we obtain $\mathcal{E}_{\sigma_0} = 1 \in  \F_1(\TP^{n+1}_{\{0\}})$.
\begin{figure}
  \centering
  \includegraphics[width=5cm,height=4.5cm]{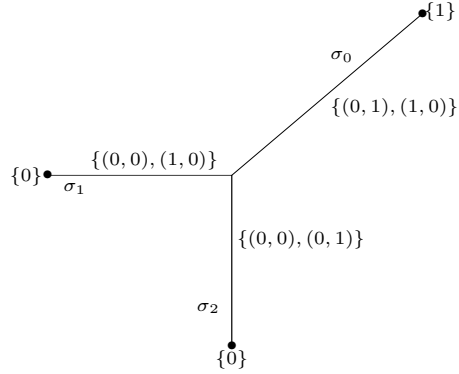}
  \put(-143,64.5){\tiny $\bullet$}
  \put(-155,64){\tiny $\left\lbrace 0\right\rbrace$}
  \put(-125,70){\tiny $\left\lbrace (0,0), (1,0) \right\rbrace$} 
   \put(-135,60){\tiny $\sigma_1$} 
  \put(-74,0){\tiny $\bullet$}
  \put(-78,-5){\tiny $\left\lbrace 0\right\rbrace$}
  \put(-70,40){\tiny $\left\lbrace (0,0), (0,1) \right\rbrace$}
    \put(-85,15){\tiny $\sigma_2$} 
  \put(-2.5,125.5){\tiny $\bullet$}
  \put(0,126){\tiny $\left\lbrace 1\right\rbrace$}
  \put(-35,90){\tiny $\left\lbrace (0,1), (1,0) \right\rbrace$}
    \put(-35,110){\tiny $\sigma_0$} 
  \caption{The extension of the real phase structure on the   tropical line from Example \ref{exrealtropline}  to faces of non-empty  sedentarity.}
\end{figure}
\end{exa} 

\begin{exa}\label{explanesed}
The real phase structure on the tropical plane  $X \subset \T P^3$ from  Example \ref{exrealtropplane} can be extended to the facets of all strata $X_I$ for $I$ a proper subset of $\{0, \dots, n+1 \}$. 
If $|I| = 1$ then $X_I$ is a tropical  line in as in Example \ref{exrealtropline}. Consider for example $I = \{1\}$ and the facet $\rho_2  = \overline{\sigma_{12}} \cap X_{\{1\}}$. 
The projection $\pi_{\sigma_{12} \rho_2}$ has kernel the first coordinate direction and therefore $\mathcal{E}_{\rho_2} = \langle e_2 \rangle \subset \F_1(\TP^3_{\{1\}})$. 
Furthermore, $X_{\{12\}}$ is a point  and $\mathcal{E}_{X_{\{12\}}} = 0$. 

\end{exa}

For any facet $\sigma$, we define the abstract vector space $\mathcal{S}_{\mathcal{E}}(\sigma)$ with generators in bijection with the elements of $\mathcal{E}_\sigma$,
$$
\mathcal{S}_{\mathcal{E}}(\sigma)=\Z_2\left\langle w_\varepsilon\mid \varepsilon\in\mathcal{E}_\sigma\right\rangle.
$$
The vector space $\mathcal{S}_{\mathcal{E}}(\sigma)$ is a linear subspace of the abstract vector space 
$
\Z_2\left\langle w_\varepsilon\mid \varepsilon\in\Z_2^{n+1-|\sed(\sigma)|}\right\rangle. 
$

\begin{definition}\label{def:signcosheaf}
Let $X$ be a non-singular real tropical hypersurface equipped with a real phase structure $\mathcal{E}$. The \emph{sign cosheaf} $\S_\mathcal{E}$ on $X$ is defined by 
\begin{equation}\label{def:S}
\S_\mathcal{E}(\tau) = \sum_{\substack{ \sigma \: \mid \: \tau \subset \sigma \\ \sed(\sigma) = \sed(\tau)}} \mathcal{S}_{\mathcal{E}}(\sigma).
\end{equation}
The maps of the cellular cosheaf 
$$i_{\sigma \tau} \colon \S_\mathcal{E}(\sigma) \to \S_\mathcal{E}(\tau)$$
are induced by natural inclusions when $\sigma$ and $\tau$ are in the same boundary stratum of $X$ and otherwise are induced by the quotients $\pi_{\sigma \tau}$ from (\ref{projmap})  composed with inclusions. 
\end{definition}

\begin{exa}\label{excosheafonP}
We describe some of the vector spaces $\mathcal{S}_{\mathcal{E}}(\tau)$ and maps between them for the  real phase structure on the   tropical plane  $X \subset \T P^3$ from  Example \ref{exrealtropplane}.

For the facets $\sigma_{01}, \sigma_{12},$ and $\sigma_{13}$ of sedentarity $0$ from Example \ref{exrealtropplane} we have, 
$$\mathcal{S}_{\mathcal{E}}(\sigma_{01}) = \langle w_{\va_3}, w_{\va_1 + \va_3}, w_{\va_1 + \va_2} , w_{\va_{2}} \rangle,$$ 
$$\mathcal{S}_{\mathcal{E}}(\sigma_{12}) = \langle w_0, w_{\va_1}, w_{\va_2} , w_{\va_{1} + \va_2} \rangle, \ \text{and} $$ 
$$\mathcal{S}_{\mathcal{E}}(\sigma_{13}) = \langle w_0, w_{\va_1}, w_{\va_3} , w_{\va_{1} + \va_3} \rangle.$$ 
Consider the one dimensional face $\tau_1$ of sedentarity $0$ and in direction $e_1$. Then we have 
$$\mathcal{S}_{\mathcal{E}}(\tau_1) = \langle w_{\varepsilon} \ | \ \varepsilon \in \mathcal{E}_{\sigma_{12}}  \cup  \mathcal{E}_{\sigma_{13}} \cup  \mathcal{E}_{\sigma_{01}} \rangle = \langle w_{0}, w_{\va_{1}},
w_{\va_{2}}, w_{\va_{1}  + \va_2}, w_{\va_{3}}, w_{\va_{1} + \va_3}  \rangle,$$
and there is an injection $i_{\sigma_{12} \tau_1} \colon \mathcal{S}_{\mathcal{E}}(\sigma_{12}) \to \mathcal{S}_{\mathcal{E}}(\tau_1)$.

For the face  $\rho_2$ from Example \ref{explanesed} we have $$\mathcal{S}_{\mathcal{E}}({\rho_2}) = \langle w_0, w_{\va_{2}} \rangle \subset \Z_2\left\langle w_\varepsilon\mid \varepsilon\in\F_1(\TP^{n+1}_{\{1\}})\right\rangle.$$
The map $i_{\sigma_{12} \rho_2} \colon \mathcal{S}_{\mathcal{E}}(\sigma_{12}) \to \mathcal{S}_{\mathcal{E}}(\rho_2)$ 

 has kernel equal to $w_{\va_{1}}$. 

\end{exa}

\begin{lemma}
Let $X$ be a non-singular real tropical hypersurface equipped with a real phase structure $\mathcal{E}$. 
If $\tau$ is a face of the stratum $X_ \rho$ of dimension $k$, the dimension of $\mathcal{S}_\mathcal{E}(\tau)$ is 
$$
\dim\mathcal{S}_\mathcal{E}(\tau)=2^{\dim X_\rho + 1}-2^{k}.
$$
\end{lemma}
\begin{proof}
It follows from the definition that $\mathcal{S}_\mathcal{E}(\tau)$ is the $\Z_2$-vector space generated by $w_{\varepsilon}$ for elements $\varepsilon$ in the set
$$
\mathcal{E}(\tau):=\bigcup_{\substack{ \sigma \: \mid \: \tau \subset \sigma \\ \text{sed}(\sigma) = \text{sed}(\tau)}} \mathcal{E}(\sigma).
$$
It follows from  \cite[Proposition 3.1]{Itenberg1997}, which is formulated in the original description of Viro's patchworking described in  Remark \ref{Viro}, 
 that $\#(\mathcal{E}(\tau))=2^{\dim X_ \rho +1}-2^{k}$. 
\end{proof}
\begin{corollary}
\label{cor:maximal}
For any face $\tau$ of $X$, we have
$$
\sum_{0\leq p\leq |\tau|} \dim \F_p(\tau)=\dim \mathcal{S}_\mathcal{E}(\tau).
$$
\end{corollary}
\begin{proof}

Let $\tau$ be a face of the stratum $X_ \rho$ of dimension $k$. It follows from Lemma \ref{lem:dimFp} that
$$
\chi_\tau(\lambda)=(1-\lambda)^{k}\left[(1-\lambda)^{\dim X_\rho +1}-(-\lambda)^{\dim X_\rho +1}\right],
$$
and then
$$
\sum_{0\leq p\leq |\tau|}\dim\F_p=\chi_{\tau}(-1)=2^{\dim X_\rho +1}-2^{k}=\dim\mathcal{S}_\mathcal{E}(\tau).
$$
\end{proof}

\begin{definition}\label{def:realhomo}
Let $X$ be a non-singular  tropical hypersurface equipped with a real phase structure $\mathcal{E}$. 
The groups of \emph{cellular $q$-chains with coefficients in $\mathcal{S}_{\mathcal{E}}$} are 
$$C_q(X ; \mathcal{S}_{\mathcal{E}}) = \bigoplus_{\dim \sigma = q} \mathcal{S}_{\mathcal{E}} (\sigma).$$
The boundary maps $\partial \colon C_q(X ; \mathcal{S}_{\mathcal{E}}) \to C_{q-1}(X ; \mathcal{S}_{\mathcal{E}}) $ are given by the direct sums of the cosheaf maps $i_{\sigma \tau}$ for $\tau\subset \sigma$. 
The \emph{real tropical  homology groups} are 
$$H_q(X;  \mathcal{S}_{\mathcal{E}}) := H_q(C_{\bullet}(X ; \mathcal{S}_{\mathcal{E}})).$$
\end{definition}

For a non-singular tropical hypersurface $X$ equipped with a real phase structure $\mathcal{E}$,  we relate the  homology of the cellular cosheaf $\S_\mathcal{E}$ to the homology of the real part $\R X_\mathcal{E}$. 

\begin{prop}\label{prop:realcosheaf}
Let $X $ be a  non-singular  tropical hypersurface equipped with a real phase structure $\mathcal{E}$.
There is an isomorphism of chain complexes 
$$C_{\bullet} (\R X_\mathcal{E}; \Z_2) \cong C_{\bullet}(X ;\S_\mathcal{E}).$$
It follows that 
$H_q (\R X_\mathcal{E}; \Z_2) \cong H_q(X ;\S_\mathcal{E})$ for all $q$. 
\end{prop}

\begin{proof}
The statement follows by comparing the cellular chain complexes $C_{\bullet}(\R X_\mathcal{E}; \Z_2)$ and $C_{\bullet}(X; \S_\mathcal{E})$. 
Firstly, we have $$C_{q}(\R X_\mathcal{E}; \Z_2) = \bigoplus_{\substack{ \tilde{\sigma} \in \R X  \\ \dim \tilde{\sigma}   =q}  }  \Z_2 \langle \tilde{\sigma} \rangle$$
and the differential 
\begin{equation}\label{eq:celldiff}
\partial \colon C_{q}(\R X_\mathcal{E}; \Z_2)  \to C_{q-1}(\R X_\mathcal{E}; \Z_2)
\end{equation}
is given componentwise by  maps
$ \tilde{\sigma} \to \sum_{ \tilde{\tau} \in \partial \tilde{\sigma}} \tilde{\tau}$. 
We can rewrite these chain groups by summing instead over the faces of $X$ 
$$C_{q}(\R X_\mathcal{E}; \Z_2)  = \bigoplus_{\substack{{\sigma} \in X  \\ \dim \sigma   =q}}\left( \bigoplus_{\substack{\tilde{\sigma} \in \R X_\mathcal{E} \\ \exists \varepsilon | \sigma^\varepsilon=\tilde{\sigma}}}  \Z_2 \langle \tilde{\sigma} \rangle \right)= \bigoplus_{\substack{{\sigma} \in X  \\ \dim \sigma   =q}} \mathcal{S}_\mathcal{E}(\sigma).$$
By Definition \ref{def:realhomo} we have  $C_q(\R X_\mathcal{E}; \Z_2) \cong C_q(X;\S_\mathcal{E})$ for all $q$. Also by the definition of the  maps $i_{\sigma \tau} \colon \mathcal{S}_{\mathcal{E}}(\sigma) \to \mathcal{S}_{\mathcal{E}}(\tau)$ for $\sigma$ and $\tau$ of dimensions $q$ and $q-1$ respectively, we see that the differentials of the chain complex $C_{\bullet}(X; \S_\mathcal{E})$ coincide with the differentials in (\ref{eq:celldiff}) above. Therefore the chain complexes are isomorphic and the isomorphism of homology groups follows. 
\end{proof}

\section{A filtration of the chain complex}
\label{section:filtration}

We begin by describing the augmentation filtration highlighted by Quillen \cite{Quillen}  on the abstract vector space $\Z_2\left\langle w_\varepsilon\mid \varepsilon\in V \right\rangle$,  where $V$ is a vector space defined over $\Z_2$. 
This same filtration  was used in \cite{Mtoric} and \cite{Hower} to give criteria for a toric variety to be maximal in the sense of the Smith-Thom inequality. 
We then adapt the filtration to when $V$ is an affine subspace and not only a vector space, and apply this to filter first the vector spaces $\S(\sigma)$ where
$\sigma$ is a top dimensional face of a tropical hypersurface. This produces filtrations of the spaces $\S(\tau)$ for any face $\tau$. 
Finally, we show that this produces a filtration of the chain complex $C_ \bullet (X; \S)$. 
Throughout this section $X$ will denote an $n$-dimensional tropical hypersurface with Newton polytope $\Delta$ contained in the tropical toric variety $\T Y_{\Delta}$. 

\subsection{The augmentation filtration}
Let $V$ be a vector space defined over $\Z_2$. 
The vector space $\Z_2\left\langle w_\varepsilon\mid \varepsilon\in V \right\rangle$ can also be considered as the group algebra $\Z_2\left[V\right]$, where the algebra structure is given by $w_\varepsilon w_\eta=w_{\varepsilon+\eta}$.
Any element of $\Z_2\left[V\right]$ can be written as $\sum a_i w_{\varepsilon_i}$, where $a_i\in\Z_2$ and $\varepsilon_i\in V$. For a subset $G \subset V$ we define $$w_G  := \sum_{\varepsilon \in G} w_{\varepsilon}.$$  

The \emph{augmentation morphism} is  given by 
$$
\begin{array}{ccccc}
\varphi & : &\Z_2\left[V\right] & \to & \Z_2 \\
 & & \sum a_i w_{\varepsilon_i} & \mapsto & \sum a_i.\\
\end{array}
$$ 
The  \emph{augmentation ideal} of $V$, denoted $\I_V$ is the  kernel of this morphism.
For all $p \geq 1$, define
$$
\I_V^p=\left\lbrace w_1\cdots w_p \mid w_1,\cdots,w_p \in\I_V \right\rbrace.
$$ Notice that  $\I_V = \I_V^1$. Since $\I_V$ is the kernel of a homomorphism it is an ideal, and for all $p \geq 1$ we have $\I_V^{p+1}\subset \I_V^p$. We obtain a  filtration of $\Z_2\left[V\right]$:
\begin{equation}
\label{augmentationfiltration}
\cdots \subset \I_V^{p} \subset \cdots \subset  \I^2_V  \subset \I_V \subset \Z_2\left[V\right].
\end{equation}
The following lemma and proposition are generalisations of \cite[Lemma 6.1 and Proposition 6.1]{Mtoric}, but we recall their proofs for  convenience.
 The Grassmannian of $p$ dimensional vector subspaces of $V$ is denoted by $\text{Gr}_p(V)$.

\begin{lemma}\label{lemma:generators}
For $p\geq 1$, the power $\I_V^p$ is additively generated by $\{w_F \ | \ F \in \Gr_p(V)\}$. 
\end{lemma} 
\begin{proof}
The proof is by induction on $p$. If $p=1$, an element $w\in \I_V$ is a sum of an even number of elements $w_{\varepsilon_1},\cdots,w_{\varepsilon_{2l}}\in\Z_2[V]$. Then we can write 
$$
w=\sum_{i|w_{\varepsilon_i}\neq w_0}w_{\left\lbrace 0,\varepsilon_i\right\rbrace},
$$
and  ${\left\lbrace 0,\varepsilon_i\right\rbrace}$ is a $1$-dimensional subspace of $V$. 

Now assume that the claim is true for $p$. Then $\I_V^{p+1}$ is additively generated by products $w_{F'} w_{\left\lbrace 0,\varepsilon\right\rbrace}$, where $F'$ is a vector subspace of $V$ of dimension $p$ and $\varepsilon\in V$. If $\varepsilon\in F'$ then $w_{F'} w_{\left\lbrace 0,\varepsilon\right\rbrace}=0$ and otherwise $w_{F'} w_{\left\lbrace 0,\varepsilon\right\rbrace}=w_F$, where 
$F=\langle F', \varepsilon \rangle$ 
is a vector subspace of $V$ of dimension $p+1$.
This completes the proof. 
\end{proof}
\begin{corollary}
If $p>\dim V$, then $\I_V^p=0$, and the filtration in  (\ref{augmentationfiltration}) is a  filtration of length $\dim V$.
\end{corollary}

\begin{proposition}\label{prop:augmentation}
For all  $p$ there is an isomorphism  $\I_V^p / \I_V^{p+1}\cong  \bigwedge^p V$.
\end{proposition}
\begin{proof}
Consider the map
$$
\begin{array}{ccccc}
f & \colon &V & \to & \Z_2\left[V\right] \\
 & & \varepsilon & \mapsto & w_0 + w_{\varepsilon} \\
\end{array}
$$
This map is not a homomorphism but $f(V)\subset \I_V$. One can consider the composition of $f$ with the quotient map $\I_V \to \I_V /\I_V^2$. We again denote this map by $f$. Then $f \colon V \to \I_V/\I_V^2$ is a homomorphism since 
$$
f(\varepsilon)+f(\eta)+f(\varepsilon+\eta)=(w_0+w_\varepsilon)(w_0+w_\eta)\in \I_V^2.
$$ 
For all $p$, define the map
$$
\begin{array}{ccccc}
f_{p} & : &V^p & \to & \I_V^p / \I_V^{p+1} \\
 & & (\varepsilon_1,\cdots,\varepsilon_p) & \mapsto & f(\varepsilon_1)\cdots f(\varepsilon_p), \\
\end{array}
$$
where $(\I_V/\I_V^2)^p$ is naturally identified with $\I_V^p/\I_V^{p+1}$.
This map is $p$-linear and alternating since 
$$
(f(\varepsilon))^2=\left[(w_0+w_\varepsilon)^2\right] = w_0^2 + w_\varepsilon^2 = w_0 + w_0 =0,
$$
and it descends to a linear map 
$$
\begin{array}{ccccc}
\hat{f}_{p} & : &\bigwedge^p V & \to & \I_V^p / \I_V^{p+1} \\
 & & \varepsilon_1\wedge\cdots\wedge\varepsilon_p & \mapsto & f(\varepsilon_1)\cdots f(\varepsilon_p). \\
\end{array}
$$
By Lemma \ref{lemma:generators}, the power $\I_V^p$ is generated by $w_F$ for all vector subspaces $F$ of $V$ of dimension $p$, so the maps $\hat{f}_p$ are surjective. 
Since $\I_V=\ker \varphi$,  its dimension is $2^{\dim V}-1$, and the set of generators $\left\lbrace w_F \mid  F \in \Gr_1(V) \right\rbrace$ is a basis of $\I_V$. Therefore, 

$$
\sum_{p=1}^{\dim V} \dim \I_V^p/\I_V^{p+1}=\dim \I_V = 2^{\dim V}-1
$$
and 
$$
\sum_{p=1}^{\dim V} \dim \bigwedge^p V= 2^{\dim V}-1,
$$
so all the map $\hat{f}_p$ are isomorphisms.
\end{proof}

\subsection{The filtration of the sign cosheaf}
Here we adapt the augmentation filtration from the last subsection to filter the vector spaces $\S(\tau)$ for $\tau$ a face of a tropical hypersurface $X$  equipped with a real phase structure. 

Let $\sigma$ be a facet of $X$. By choosing a vector $\theta$ in  the affine hyperplane $\mathcal{E}_\sigma$, we obtain an identification
$$
\psi_\theta \colon \S_\mathcal{E}(\sigma) \cong \Z_2[\overrightarrow{\mathcal{E}_\sigma}],
$$
where $\overrightarrow{\mathcal{E}_\sigma}$ denotes the vector subspace of $\Z_2^{n+1}$ parallel to $\mathcal{E}_ \sigma$. Transporting the augmentation filtration of $\Z_2[\overrightarrow{\mathcal{E}_\sigma}]$ by the isomorphism $\psi_\theta$, one obtains a filtration of $\S_\mathcal{E}(\sigma)$
$$
0=\mathcal{K}_{n+1}(\sigma)\subset\mathcal{K}_{n}(\sigma)\subset\cdots\subset\mathcal{K}_{0}(\sigma)=\S_\mathcal{E}(\sigma).
$$
In the following lemma, we show that this filtration does not depend on the choice of an element $\theta$ in $\S_\mathcal{E}(\sigma)$ we choose. Let $\Aff_p(\mathcal{E}_\sigma)$ denote the space of all $p$-dimensional affine subspaces of  $\mathcal{E}_{\sigma}$.
\begin{lemma}
For any facet $\sigma$ of $X$ one has 
$$\mathcal{K}_p(\sigma) =  \langle  w_{G} \ | \ G \in \Aff_p(\mathcal{E}_\sigma)   \rangle $$
\end{lemma} 
\begin{proof}
Recall that we choose an element $\theta$ of the affine hyperplane $\mathcal{E}_\sigma$. By Lemma \ref{lemma:generators}, the vector space $\mathcal{K}_p(\sigma)$ is generated by the $w_{G}$ for all affine subspaces $G$ of $\mathcal{E}_\sigma$ of dimension $p$ passing through $\theta$. Let $G$ be an affine subspace of $\mathcal{E}_\sigma$ of dimension $p$ not passing through $\theta$, and let $H$ be any affine hyperplane $H$ of $G$. Since we are over $\Z_2$, one has $G=H\cup H'$, where $H'$ is the affine hyperplane of $G$ parallel to $H$. Denote by $H_1$ the affine subspace of $\mathcal{E}_\sigma$ parallel to $H$ and passing through $\theta$. Then one has
$$
w_G=w_{H\cup H_1}+w_{H'\cup H_1},
$$ 
and the lemma is proved.
\end{proof}

\begin{definition}
\label{def:Kp}
Let $X$  be a real tropical hypersurface with real phase structure $\mathcal{E}$. For all $p$, we  define a collection of cosheaves $\K_p$ on $X$. 
For  $\tau$  a face of $X$ of sedentarity $I$, let  
$$\K_p(\tau) = \sum_{\sigma \supset \tau} \mathcal{K}_p(\sigma) \subset \S_\mathcal{E}(\tau)$$ 
where the sum is over facets $\sigma$ of $X_I$. 
The cosheaf maps $\K_{p}(\tau_1) \to \K_{p}(\tau_2)$ for $\tau_2 \subset \tau_1 $ are the restrictions of the maps $i_{\tau_1 \tau_2} \colon \S_{\mathcal{E}}
(\tau_1) \to \S_{\mathcal{E}}(\tau_2). $
\end{definition}

For each face $\tau$ of $X$ we obtain a filtration of $\S_\mathcal{E}(\tau)$ given by
\begin{equation}\label{filttau} \K_{n}(\tau) \subset \dots \subset  \K_{2}(\tau) \subset \K_{1}(\tau) \subset \K_{0}(\tau) =\S_\mathcal{E}(\tau) .
\end{equation}
If  $\tau_2 \subset \tau_1$, then  the facets adjacent to $\tau_1$ are a subset of the facets adjacent to $\tau_2$.  It follows that $i_{\tau_1 \tau_2}(\K_p(\tau_1) ) \subset \K_p(\tau_2)$ 
so that the cosheaf maps for $\K_p$ are well-defined.

\begin{exa}\label{exK1}
Any two  vectors $w_{\va_1}, w_{\va_2}$ for $\va_1, \va_2 \in \F_1(\TP^{n+1}_I)$ are on an affine line. Moreover, these are the only points over $\Z_2$ contained on the line. Therefore, every facet $\sigma$ of sedentarity $I$ of a real tropical hypersurface $X$, the 
 vector subspace $\K_1(\sigma) \subset \mathcal{S}_\mathcal{E}(\sigma)$ is generated by $w_{\va_1}+ w_{\va_2}$, for any  vectors $\va_1, \va_2 \in \mathcal{E}(\sigma)$. This implies that  $\K_1(\sigma)$ is the hyperplane inside $\mathcal{S}_\mathcal{E}(\sigma)$ defined by the linear form 
 $\sum_{\va \in \mathcal{E}_\sigma   } x_{\va}  = 0$ where the $x_{\va}$'s form a  dual basis to the $w_{\va}$'s. 

For a face $\tau$ of higher codimension, the space $\K_1(\tau)$ is also an hyperplane inside $\mathcal{S}_\mathcal{E}(\tau)$ defined by the linear form 
 $\sum_{\va \in \mathcal{E}_\tau  } x_{\va}  = 0$. By definition of $\K_1(\tau)$ we have 
that 
$ \K_1(\tau)$ is contained in the hyperplane defined by $ \sum_{\va \in \mathcal{E}_\tau  } x_{\va}  = 0$
To prove the reverse inclusion, it is enough to show that $w_{\va}+w_{\va'}\in \K_1(\tau)$, for any $w_{\va},w_{\va'}\in \mathcal{S}_\mathcal{E}(\tau)$. Let  $\sigma$ and $\sigma'$ be two facets of $X$ containing $\tau$ such that $\va\in\mathcal{E}_{\sigma}$ and $\va ' \in\mathcal{E}_{\sigma'}$. The intersection $\sigma\cap\sigma'$ is a face of codimension either one or two. 
 If it is a face of codimension one, then by the  condition on a real tropical  structure in Definition \ref{realphase}, there exists $\va_1\in\mathcal{E}_\sigma\cap\mathcal{E}_{\sigma'}$. But then $w_{\va}+w_{\va'}=(w_{\va}+w_{\va_1})+(w_{\va_1}+w_{\va_2})\in\K_1(\tau)$. If $\sigma\cap\sigma'$ is a face of codimension two, then there exists a facet $\sigma''$ such that $\sigma\cap\sigma''$ and $\sigma'\cap\sigma''$ are of codimension $1$ and $\sigma\cap\sigma'\cap\sigma''=\sigma\cap\sigma'$. Then, there exist $\va_1\in\mathcal{E}_\sigma\cap\mathcal{E}_{\sigma''}$ and $\va_2\in\mathcal{E}_{\sigma''}\cap\mathcal{E}_{\sigma'}$ such that 
$$ 
w_{\va}+w_{\va'}=(w_{\va}+w_{\va_1})+(w_{\va_1}+w_{\va_2})+(w_{\va_2}+w_{\va'})\in\K_1(\tau).
$$
This shows that $\K_1(\tau)$ is also a hyperplane inside $\S_{\mathcal{E}}(\tau)$ for all faces $\tau$. 
\end{exa}

\begin{exa}\label{exK2plane}
For the real tropical plane from Example \ref{exrealtropplane} we describe the filtration in $(\ref{filttau})$ for some faces. 
Following Example \ref{exK1}, for every facet $\sigma_{ij}$ of $X$ the vector space $\K_1(\sigma_{ij})$ is  of codimension one in $\mathcal{S}_{\mathcal{E}}(\tau)$. 
For any facets  $\sigma_{ij}$ of $X$ the vector space $\mathcal{S}_{\mathcal{E}}(\sigma)$ is two dimensional. Therefore, the only element 
in $\Aff_2(\mathcal{S}_{\mathcal{E}}(\sigma))$ is the whole vector space itself. This implies that $\K_2(\sigma_{ij}) = \langle w_{\mathcal{S}_{\mathcal{E}}(\sigma_{ij} )} \rangle $, in particular it is one dimensional. 
For instance for $\sigma_{12}$ we have,
$$\mathcal{S}_{\mathcal{E}}(\sigma_{12}) = \langle w_0, w_{\va_1}, w_{\va_2}, w_{\va_1 + \va_2} \rangle, $$
 $$\K_{1}(\sigma_{12}) = \langle w_0 +  w_{\va_1}, w_0 +  w_{\va_2}, w_0 +  w_{\va_1 + \va_2} \rangle,  \text{and} $$
$$  \K_{2}(\sigma_{12}) = \langle w_0 +  w_{\va_1}  +  w_{\va_2} +   w_{\va_1 + \va_2} \rangle.$$

For the face $\tau_1$ from Example \ref{excosheafonP}, since $\K_{1}(\tau_{1})$ is generated by $\K_1(\sigma_{01}),$ $\K_1(\sigma_{12}),$ and $\K_1(\sigma_{13})$, we have 
$$\K_{1}(\tau_{1}) = \langle w_0 +  w_{\va_1}, w_0 +  w_{\va_2}, w_0 +  w_{\va_1 + \va_2},  w_0 +  w_{\va_3}, w_0 +  w_{\va_1 + \va_3} 
\rangle.$$
For $p= 2$ we have 
$$\K_{2}(\tau_{1}) = \langle w_{\mathcal{S}_{\mathcal{E}}(\sigma_{01})} ,  w_{\mathcal{S}_{\mathcal{E}}(\sigma_{12})} 
\rangle,$$
since $w_{\mathcal{S}_{\mathcal{E}}(\sigma_{01})} + w_{\mathcal{S}_{\mathcal{E}}(\sigma_{12})} + w_{\mathcal{S}_{\mathcal{E}}(\sigma_{13})} = 0$. 
\end{exa}

\begin{lemma} \label{lemma:isoK_p}
For any face $\tau$ of $X$, there is an isomorphism 
$$
\F_p(\tau) \cong  \K_p(\tau)/\K_{p+1}(\tau).
$$
\end{lemma}
\begin{proof}
Extend the map $\widehat{f}_p$ from to the proof of Proposition \ref{prop:augmentation} to a map defined on any face $\tau$ of $X$:
$$
\widehat{f}_p(\tau):\F_p(\tau)\longrightarrow \K_p(\tau)/\K_{p+1}(\tau).
$$ This map is again surjective and it follows from Corollary \ref{cor:maximal} that it is a isomorphism.
\end{proof}

From the isomorphism from Lemma \ref{lemma:isoK_p}  homomorphisms $\K_p(\tau) \to \F_p(\tau)$, which we call the Viro homomorphisms following \cite{Degtyarev}. 
\begin{definition}\label{def:bvp}
For any face $\tau$ of $X$, define the Viro homomorphisms $bv_p \colon \K_p(\tau) \to \F_p(\tau)$ as the composition of the quotient map $$\K_p(\tau)\rightarrow \K_p(\tau) / \K_{p+1}(\tau)$$ with the inverse of the isomorphism $\widehat{f}_p(\tau)$.
\end{definition}

If $\tau$ is a face  of $X$, then $\K_p(\tau)$ is also generated by vectors of the form $w_{G}$, where $G$ is an element of $\text{Aff}_p(\mathcal{E}_{\sigma})$ for some top dimensional face $\sigma$ containing $\tau$.
The Viro map on the generators is $$bv_p(w_{G}) = v_1 \wedge \dots \wedge v_p,$$ where
$v_1, \dots, v_p$ is a basis of the vector space parallel to the affine space $G$.

\begin{proposition} \label{prop:exactcosheaf}
For all faces $\tau \subset \sigma$ of $X$, the following diagram is commutative
\begin{align}\label{diagcomm}
\begin{xy}
\xymatrix{
 0 \ar[r] &  \K_{p+1}(\sigma) \ar[rr]^{i} \ar[d]^{i_{\sigma \tau}} && \K_{p}(\sigma) \ar[rr]^{bv_p}  	 \ar[d]^{i_{\sigma \tau}}			&& \F_{p}(\sigma) \ar[d]^{i_{ \sigma \tau}}  \ar[r] & 0  \\
0 \ar[r] & \K_{p+1}(\tau) \ar[rr]^{i}&& \K_{p}(\tau) \ar[rr]^{bv_p}           && \F_{p}(\tau)   \ar[r] & 0.
}
\end{xy}
\end{align}
\end{proposition}

\begin{proof}
The exactness of  the rows follows from Lemma \ref{lemma:isoK_p}.
Since the augmentation morphism commutes with linear projections and inclusions, the left-hand square is commutative. The commutativity of the square on the right follows from the description of $bv_p$ on the generators. 
\end{proof}

The cellular $q$-chains with coefficients in $\K_p$ are defined by 
$$C_{q}(X ; \K_p)   = \bigoplus_{\dim \sigma = q} \K_p(\sigma).$$ 
Thanks to the  commutativity of the left hand square of the  diagram in Proposition \ref{prop:exactcosheaf}, there is the    complex of relative chains $$C_{\bullet}(X ; \K_{p}, \K_{p+1}):= C_{\bullet}(X; \K_{p})/C_{\bullet}(X; \K_{p+1}).$$
We let $H_q(X ; \K_{p}, \K_{p+1}) $ denote the $q$-th homology group of this complex. 

\begin{cor}\label{cor:relativenadtrop}
For all $p$ and $q$ we have  isomorphisms
$$H_q(X ; \K_{p}, \K_{p+1})  \cong H_q(X ;\F_p).$$
\end{cor}

\begin{proof}
There is an  isomorphism $bv_{p} \colon C_{q}(X; \K_p, \K_{p+1} ) \to  C_{q}(X ; \F_{p})$ for each $q$. The commutativity on the right hand side of Proposition \ref{prop:exactcosheaf} implies that  $bv_{p} $ induces an isomorphism of complexes $C_{\bullet}(X; \K_p, \K_{p+1} ) \to  C_{\bullet}(X ; \F_{p})$. Since the complexes are isomorphic, so are their homology groups and this proves the  statement of the corollary. 
\end{proof}

\begin{prop}
\label{prop:firstpage}
The first page of the spectral sequence associated to the filtration of the chain complex  $C_{\bullet}(X ;\S_{\mathcal{E}})$ by  the chain complexes $C_{\bullet}(X ;\K_p)$ has  terms 
$$E^1_{q,p} \cong  H_q(X; \F_p).$$
\end{prop}

\begin{proof}
Proposition \ref{prop:exactcosheaf} implies that the chain complexes $C_{\bullet}(X ; \K_p)$ filter the  chain complex $C_{\bullet}(X ; \mathcal{S}_{\mathcal{E}})$ from Definition \ref{def:realhomo}
 $$0 \subset C_{\bullet}(X ; \K_n)  \subset \dots \subset C_{\bullet}(X ; \K_1) \subset C_{\bullet}(X ; \mathcal{S}_{\mathcal{E}}).$$
This is a finite filtration of a complex of finite dimensional vector spaces, therefore  the spectral sequence associated to this filtration converges \cite[Theorem 2.6]{McCleary}. 
By definition, the  first page of the spectral sequence of  the filtered complex consists of the relative chain groups,
$$E^1_{q,p} \cong H_{q}(X ;\K_p, \K_{p+1}).$$ 
Then the proposition follows from Corollary  \ref{cor:relativenadtrop}. 
\end{proof}
\begin{proof}[Proof of Theorem \ref{thm:toricvar}]

The pages of a  spectral sequence satisfy  $\dim E^{\infty}_{q,p} \leq \dim E^r_{q,p}$ for all $r$.  By Propositions \ref{prop:realcosheaf},   \ref{prop:firstpage}, and the convergence of the spectral sequence associated to the filtration  we obtain 
$$\dim H_q (\R V) = \dim H_q(X; \S_{\mathcal{E}})  = \sum_{p = 0}^{\dim X}  E^{\infty}_{q,p}  \leq  \sum_{p = 0}^{\dim X}  E^{1}_{ q,p} =  \sum_{p = 0}^{\dim X} H_q(X; \mathcal{F}_p).$$

When the tropical hypersurface $X$ is contained in a partial compactification of the torus  $\T Y_{\Delta}^o$ corresponding to a subfan of the dual fan of the Newton polytope  $\Delta$
the filtration of the chain complex $C_{\bullet}(\overline{X} ;\S_{\mathcal{E}})$ can be restricted to the cells contained in $X= \overline{X} \cap Y_{\Delta}^o$ to give a filtration of $C_{\bullet}^{BM}(X ;\S_{\mathcal{E}})$. Variants of Propositions \ref{prop:realcosheaf} and  \ref{prop:firstpage} also hold in the non-compact case and the argument given above completes the proof. 
\end{proof}

\section{Proof of Theorem \ref{thm:main}}
\label{section:mainproof}
\begin{lemma}\label{lemvanishKS}
Let $X$ be a compact tropical non-singular  hypersurface of dimension $n$ in the non-singular tropical toric variety $\T Y_\Delta$, where $\Delta $ is the Newton polytope of $X$. 
Then $$\dim H_q (X ; \F_p)  = 0 $$ unless $p + q = n$ or $p = q$. 
\end{lemma}

\begin{proof}
Let $V$ be a complex non-singular hypersurface of the same dimension and Newton polytope as $X$ considered in the complex toric variety $Y_{\Delta}$. 
The Lefschetz Hyperplane Section  Theorem together with Poincar\'e duality for $\C V$ implies that $h^{p, q}(\C V) = 0$ unless  $p + q = n$ or $p = q$. 
The statement of the lemma now follows by applying Theorem \ref{thm:hodgemod2}.
\end{proof}

\begin{proof}[Proof of Theorem \ref{thm:main}]
It follows from the statement of Theorem \ref{thm:toricvar} that
$$b_q(\R V)  \leq  \sum_{p = 0}^d \dim H_q(X; \F_p).$$
By Lemma \ref{lemvanishKS}, the sum on the right hand side 
is equal to $\sum_{p = 0} h^{p, q}(\C V)$ which completes the proof.
\end{proof}

\section{Going further in the  spectral sequence}
\label{sec:spectralsequence}

In addition to bounding the Betti numbers of real hypersurfaces close to a non-singular tropical limit, the spectral sequence provides immediate criteria for the optimality of the  bounds on individual Betti numbers from Theorem \ref{thm:main}, in addition to the criterion for maximality in the sense of the  Smith-Thom inequality from Theorem \ref{thmmax}.

\begin{thm}\label{thmindivid}
Let $V $  be  a compact  hypersurface with non-singular  Newton polytope near a non-singular tropical limit, then  
the $q$-th Betti number of $\R V$ attains the bound in Theorem \ref{thm:main} if and only if all of the following maps are zero
\begin{enumerate}
\item   when  $q=n/2$,
$$\partial_1 \colon E^1_{q, q} \to E^1_{q-1, q+1} \qquad \text{and} \qquad \partial_1 \colon E^1_{q+1, q-1}  \to E^1_{q, q},$$

\item  when  $q < n/2$,  
$$\partial_1 \colon E^1_{q, n-q} \to E^1_{q-1, n-q+1}, \qquad  \partial_1 \colon E^1_{q+1, n-q-1}  \to E^1_{q, n-q}, \qquad \text{and} $$
$$\partial_{2q  -n +  1} \colon  E^{2q  -n  + 1}_{q+1, q-r} \to E^{2q  -n  + 1}_{q, q}.$$
\end{enumerate}
\end{thm}

\begin{remark}  
If $V $  is a compact  hypersurface in a non-singular toric variety   near a non-singular tropical limit, then  the real point set $\R V$ is a smooth $n$-dimensional  manifold and its Betti numbers over $\Z_2$ satisfy Poincar\'e duality. 
This ensures that $$b_i(\R V) = b_{n-i}(\R V).$$ Therefore in order to determine all of the  Betti numbers of $\R V$  we only need to determine the Betti numbers $b_q(\R V) $ for $q \leq n/2$. 
\end{remark}

\begin{lemma}\label{lemvanishingmaps}
Let $V $  be  a compact  hypersurface with non-singular  Newton polytope near a non-singular tropical limit.
 The only possible non-zero differentials of the spectral sequence $(E^{\bullet}_{\bullet, \bullet}, \partial^{\bullet})$ 
 are  \begin{align}
  \label {map1} \partial_1 \colon E^1_{q, p} \to E^1_{q-1, p+1} & \quad  \text{ for }\quad p + q = n \\
 \label{maprn} \partial_r \colon E^r_{q+1, q-r} \to E^r_{q, q}& \quad  \text{ for } \quad  r = 2q -n +  1 \\
\label{maprqq} \partial_r \colon E^r_{q, q} \to E^r_{q-1, q+r} & \quad \text{ for }  \quad r  = n -2q +1.
\end{align}
\end{lemma}

\begin{proof}
If  a boundary  map $\partial_r \colon E^r_{q, p} \to E^r_{q-1, p+r}$ is non-zero, then necessarily  both  $E^r_{q, p}$ and  $ E^r_{q-1, p+r} $   must be non-zero. 
This implies that both $E^1_{q, p} \cong H_q(X ; \F_p)$ and  $ E^1_{q-1, p+r}  \cong  H_{q-1}(X ; \F_{p+r})$   must be non-zero.  But Lemma \ref{lemvanishKS} implies that 
$E^1_{q, p} \cong H_q(X ; \F_p) = 0$ unless $p + q = n$ or $p = q$.

Case 1: Suppose $p + q = n$, then for  $E^r_{q-1, p+r}$ to be non-zero we must have either $q-1 +  p+r = n$ or $q-1 = p+r$. 
In the first case $r =1$. In the second case $r = 2q - n- 1$. These are the maps listed in (\ref{map1}) and (\ref{maprn}) in the statement of the lemma. 

Case 2: If $p = q$, then we have the non-zero map in (\ref{map1}) when $r=1$. If $r >1$, then $q-1 + q+r = n$ so that we find the condition in (\ref{maprqq}) above. 
 This completes the proof of the lemma. 
\end{proof}

\begin{proof}[Proof of Theorem \ref{thmindivid}]
The $q$-th Betti number attains the bounds in Theorem \ref{thm:main} if and only if for all $p$ the maps  $$\partial_r \colon E^r_{q, p} \to E^r_{q-1, p+r} \quad \text{and} \quad  
\partial_r \colon E^r_{q+1, p-r} \to E^r_{q, p},$$ are zero for all $r$. 
The theorem follows from the list of possible non-zero differential maps in Lemma \ref{lemvanishingmaps}. 
\end{proof}

\begin{example}
Applying Theorem  \ref{thmindivid} to the case $n=2$ implies that for all $r$ the only non-zero differentials of the spectral sequence are 
\begin{equation}\label{ex:surfacedifferentials}
\partial_1 \colon H_2(X ; \K_0, \K_1) \to H_1(X ; \K_1, \K_2)  \quad \text{and} \quad 
 \partial_1 \colon H_1(X ; \K_1, \K_2) \to H_0(X ; \K_2).
\end{equation}
Recall that Corollary \ref{cor:sign} relates the signature of $\C V$ to the Euler characteristic to $\R V$ for a non-singular real hypersurface obtained from a primitive patchworking. Combining this with Poincar\'e duality for $\R V$, Serre duality for $\C V$, and also the Lefschetz Hyperplane Theorem for $\C V$ and $\C Y_{\Delta}$ we obtain the following equality, 
$$2b_0(\R V) -b_1(\R V)  =  2 + 2h^{2,0}(\C V) - h^{1,1}(\C V).$$ 
  Therefore a compact surface in  a  three dimensional toric variety obtained by primitive patchworking is 
maximal if and only if one of the maps in (\ref{ex:surfacedifferentials}) is zero. 

\end{example}

\begin{example}
For $n=4$ we show the first pages of the spectral sequence. The first page on the left below has terms $E^1_{q,p} \cong H_q(X;\F_p)$. 
$$
  {\tiny
    \xymatrix @!0 @R=6mm @C=.6cm {
    &   &   &   & \Z_2 &   &   &   & \\
    &   &   & 0 &   & 0 &   &   & \\
    &   & 0 &   & \Z_2 \ar[ll] &   & 0 \ar[ll] &   & \\
    & 0 &   & 0 &   & 0 &   & 0 & \\
  E^1_{0,4}  &  & E^1_{1,3} \ar[ll] & & E^1_{2,2} \ar[ll] &  &E^1_{3,1} \ar[ll] & & E^1_{4,0} \ar[ll] & \\
  & 0 &   & 0 &   & 0 &   & 0 & \\
  &   & 0 &   & \Z_2 \ar[ll] &   & 0  \ar[ll]&   & \\
  &   &   & 0 &   & 0 &   &   & \\
  &   &   &   & \Z_2 &   &   &   & \\
  }}
{\tiny
  \xymatrix @!0 @R=6mm @C=.6cm {
    &   &   &   & \Z_2 &   &   &   & \\
    &   &   & 0 &   & 0 &   &   & \\
    &   & 0 &   & \Z_2 \ar[llld] &   & 0 &   & \\
    & 0 &   & 0 &   & 0 &   & 0 & \\
 E^{2}_{0,4}  &  & E^{2}_{1,3} & & E^2_{2,2} \ar[llld] &  & E^2_{3,1} \ar[llld] & & E^2_{4,0} \ar[llld] & \\
  & 0 &   & 0 &   & 0 &   & 0 \ar[llld] & \\
  &   & 0 &   & \Z_2 &   & 0 &   & \\
  &   &   & 0 &   & 0 &   &   & \\
  &   &   &   & \Z_2 &   &   &   & \\
  }}
  $$
Notice that all differentials on the second page are trivial since the conditions (\ref{maprn}) and (\ref{maprqq})
 in Lemma \ref{lemvanishingmaps} cannot be satisfied for $n=4$ and $r=2$. Therefore   $E^3_{q,p} =  E^2_{q,p}$ and the arrows of the third page are depicted on the left. 
\\
$$
{\tiny  \xymatrix @!0 @R=6mm @C= 0.6cm {
    &   &   &   & \Z_2 &   &   &   & \\
    &   &   & 0 &   & 0 &   &   & \\
    &   & 0 &   & \Z_2 \ar[lllldd] &   & 0 \ar[lllldd] &   & \\
    & 0 &   & 0 &   & 0 &   & 0 & \\
 E^{3}_{0,4}  &  & E^{3}_{1,3} & & E^3_{2,2} &  & E^3_{3,1} \ar[lllldd] & & E^3_{4,0} \ar[lllldd] & \\
  & 0 &   & 0 &   & 0 &   & 0  & \\
  &   & 0 &   & \Z_2 &   & 0 &   & \\
  &   &   & 0 &   & 0 &   &   & \\
  &   &   &   & \Z_2 &   &   &   & \\
  }}
  {\tiny
  \xymatrix @!0 @R=6mm @C=0.6cm {
    &   &   &   & \Z_2 &   &   &   & \\
    &   &   & 0 &   & 0 \ar[lllllddd] &   &   & \\
    &   & 0 &   & E^4_{1,1}  &   & 0  &   & \\
    & 0 &   & 0 &   & 0 &   & 0 & \\
 E^{4}_{0,4}  &  & E^{4}_{1,3} & & E^4_{2,2} &  & E^4_{3,1} & & E^4_{4,0} \ar[lllllddd] & \\
  & 0 &   & 0 &   & 0 &   & 0  & \\
  &   & 0 &   & E^4_{3,3} &   & 0 &   & \\
  &   &   & 0 &   & 0 &   &   & \\
  &   &   &   & \Z_2 &   &   &   & \\
  }}
$$
On the right hand side above is the fourth page of the spectral sequence. Here all differentials are zero, moreover for $r \geq 4$ all differentials are zero by Lemma \ref{lemvanishingmaps}. 

\comment{
The third page:
$$
  \xymatrix @!0 @R=8mm @C=1cm {
    &   &   &   & \Z_2 &   &   &   & \\
    &   &   & 0 &   & 0 &   &   & \\
    &   & 0 &   & \Z_2 \ar[lllldd]^{\partial_3} &   & 0 \ar[lllldd]^{\partial_3} &   & \\
    & 0 &   & 0 &   & 0 &   & 0 & \\
 E^{3}_{0,4}  &  & E^{3}_{1,3} & & E^3_{2,2} &  & E^3_{3,1} \ar[lllldd]^{\partial_3} & & E^3_{4,0} \ar[lllldd]^{\partial_3} & \\
  & 0 &   & 0 &   & 0 &   & 0  & \\
  &   & 0 &   & \Z_2 &   & 0 &   & \\
  &   &   & 0 &   & 0 &   &   & \\
  &   &   &   & \Z_2 &   &   &   & \\
  }.
$$
The fourth page:
$$
  \xymatrix @!0 @R=8mm @C=1cm {
    &   &   &   & \Z_2 &   &   &   & \\
    &   &   & 0 &   & 0 \ar[lllllddd]^{\partial_4} &   &   & \\
    &   & 0 &   & E^4_{1,1}  &   & 0  &   & \\
    & 0 &   & 0 &   & 0 &   & 0 & \\
 E^{4}_{0,4}  &  & E^{4}_{1,3} & & E^4_{2,2} &  & E^4_{3,1} & & E^4_{4,0} \ar[lllllddd]^{\partial_4} & \\
  & 0 &   & 0 &   & 0 &   & 0  & \\
  &   & 0 &   & E^4_{3,3} &   & 0 &   & \\
  &   &   & 0 &   & 0 &   &   & \\
  &   &   &   & \Z_2 &   &   &   & \\
  }.
$$
All differentials are also $0$ here.
}
\end{example}

\section{Case of plane curves}
\label{sec:Haas}

In this section we explicitly describe the only possibly non-zero differential map in the spectral sequence in the case of curves. In this case, Viro's  primitive patchworking construction, equivalently, the real phase structures on tropical curves from Section  \ref{realphasesub}, can be reformulated in terms of  \emph{admissible twists}.  

Given a compact  non-singular tropical curve $C$ in a tropical toric surface there is another equivalent way of describing a real phase structure on $C$ in terms choosing a subset of \emph{twisted edges} of the bounded edges of $C \cap \R^2$  satisfying  an admissibility condition.
A collection $T$ of bounded edges of a tropical curve $C \cap \R^2$ is \emph{admissible} if for all $\gamma \in H_1(C; \mathcal{F}_0)$ we have 
$$\sum_{e \in T \cap \text{Supp}(\gamma)} v_e = 0 \in \Z^2_2,$$
where $v_e$ is the primitive integer direction of the edge $e$. 
The edges of $T$ are called twisted edges because of how the real algebraic curve $\R V$ near the tropical limit $C$ behaves under the logarithm map.  
See the right hand-side of Figure \ref{fig:realcubic2}.

Let  $C$ be a non-singular compact tropical curve  with a real phase structure $\mathcal{E}$. For   a bounded edge $e$ of $C$,  its  symmetric copy, $e^{\varepsilon}$ in $\R C_{\mathcal{E}}$,  is adjacent to two other edges $e_1^{\varepsilon}$, $e_2^{\varepsilon}$ of $\R C_{\mathcal{E}}$ which are also contained in the quadrant corresponding to $\varepsilon$. The twisted edges for a  real phase structure $\mathcal{E}$ correspond to those edges $e$ of $C$ for which $e_1$, $e_2$ are not contained in a closed half space  of $\R^2(\varepsilon)$ whose boundary contains $e^{\varepsilon}$.
A detailed description of this approach can be found in \cite[Section  3.2]{BIMS}.

Using the  twist formulation  we describe  explicitly the map 
$$
\partial_1:H_1(C;\F_0)\rightarrow H_0(C;\F_1)
$$
arising from  the spectral sequence on the chain level when the curve $C$ is compact.  In this case, both  of the above homology groups are isomorphic to $\Z_2^g$, where $g$ is the first Betti number of  $C$. 

\begin{example}\label{ex:realcubic}
Figure \ref{fig:realcubic}, shows a  non-singular plane tropical cubic with a twist-admissible set of edges, and the image by coordinatewise logarithm map $\mathrm{Log}$ of the real part $\R V$ of the curve  $V$  which is defined by the polynomial ${\bf P}_t$ from (\ref{Viropoly}) for $t$ sufficiently large. 
Figure \ref{fig:realcubic2}  depicts $\R C_{T}$. Notice that this curve   is maximal in the sense of Harnack's inequality, namely $b_0 (\R V) = g(\C V) +1$. 
\begin{figure}
 \begin{minipage}[l]{.46\linewidth}
  \centering
  \includegraphics[width=5cm,height=4.5cm]{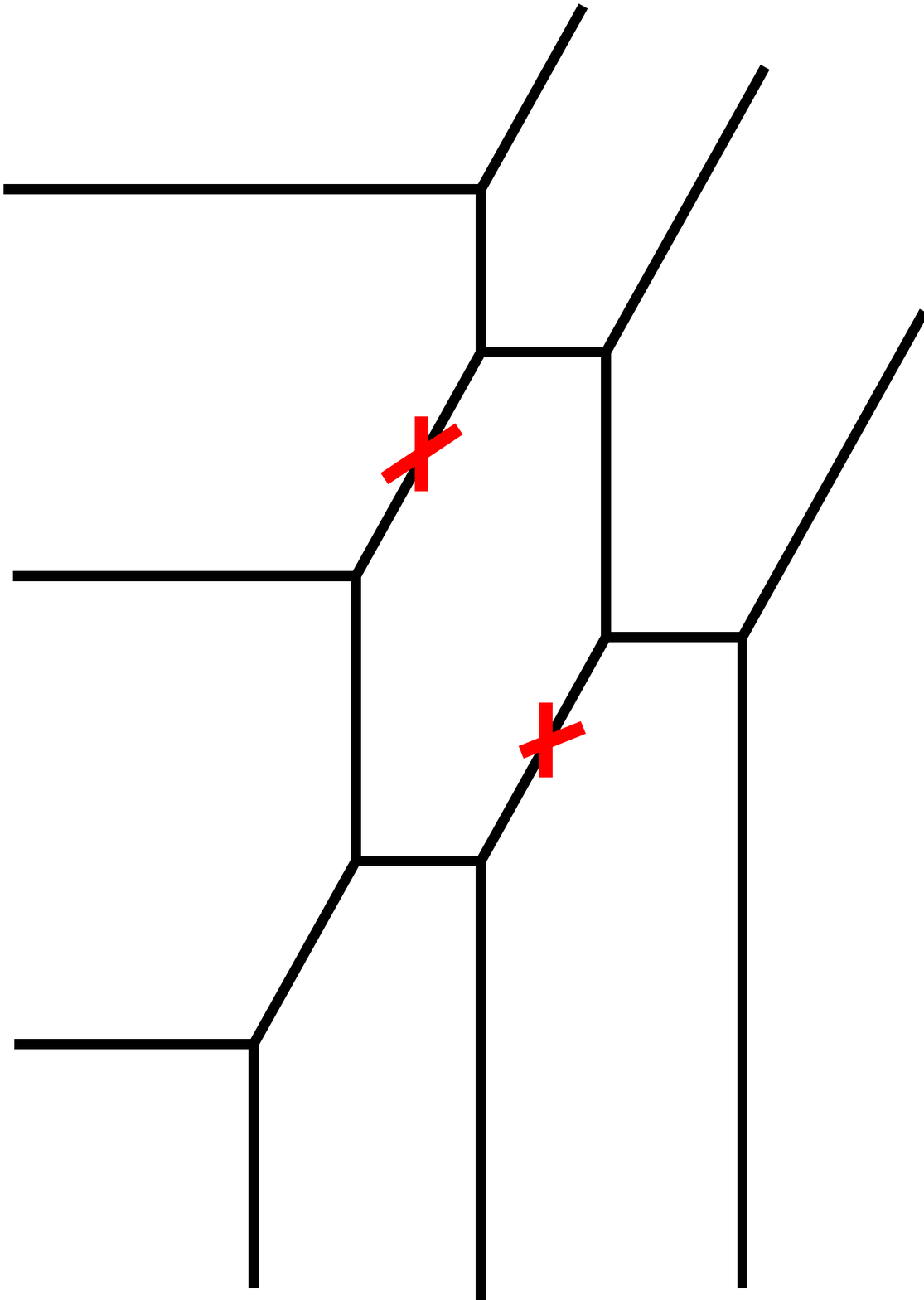}
 \end{minipage} \hfill
 \begin{minipage}[l]{.46\linewidth}
  \centering \includegraphics[width=5cm,height=5cm]{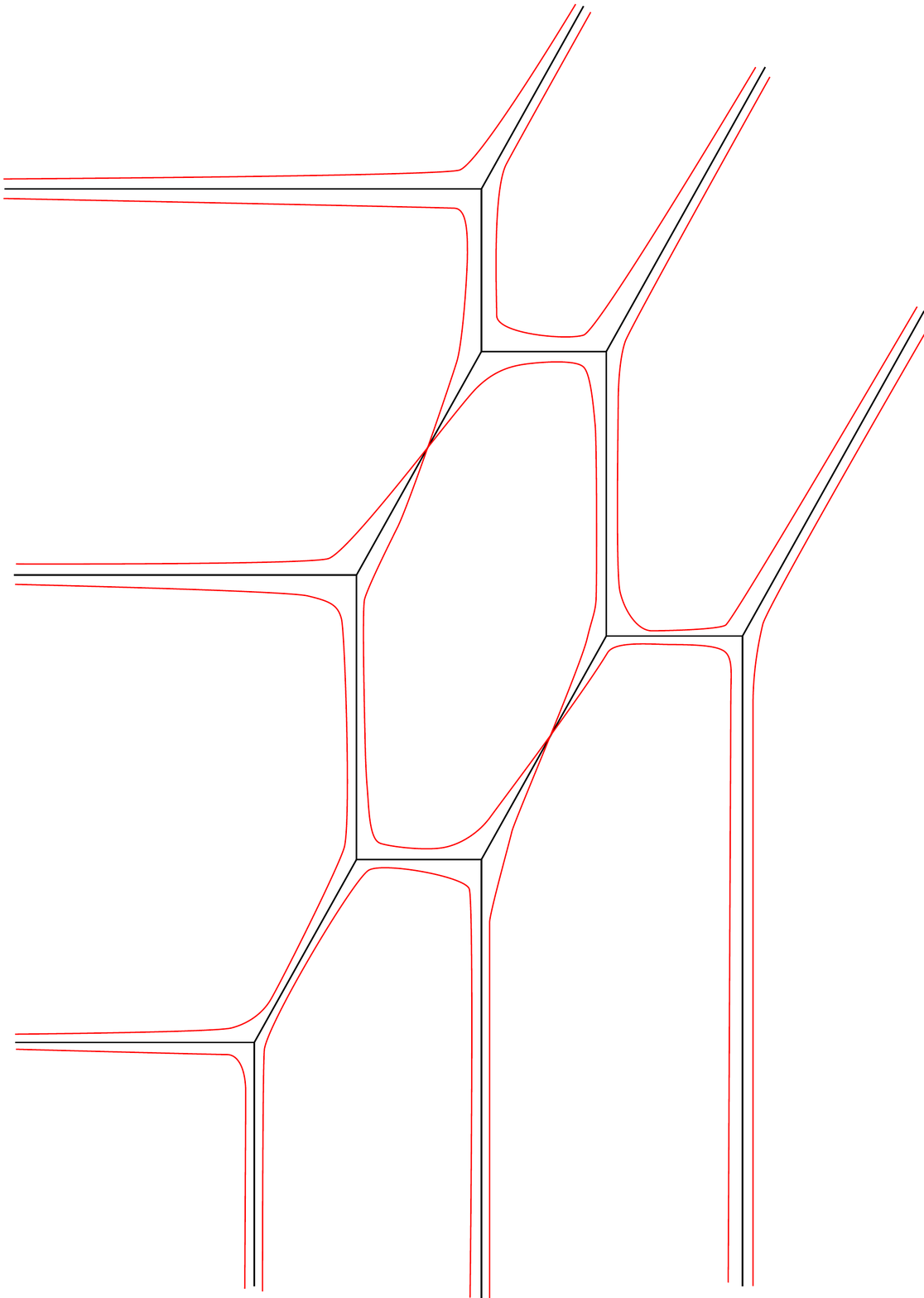}
 \end{minipage}
  \caption{On the left is a non-singular cubic with a twist-admissible set of edges. On the right hand side is the image by the coordinatewise logarithm map of $\R V$. \label{fig:realcubic}} 
\end{figure}
\begin{figure}
\centering
  \includegraphics[width=5cm,height=4.5cm]{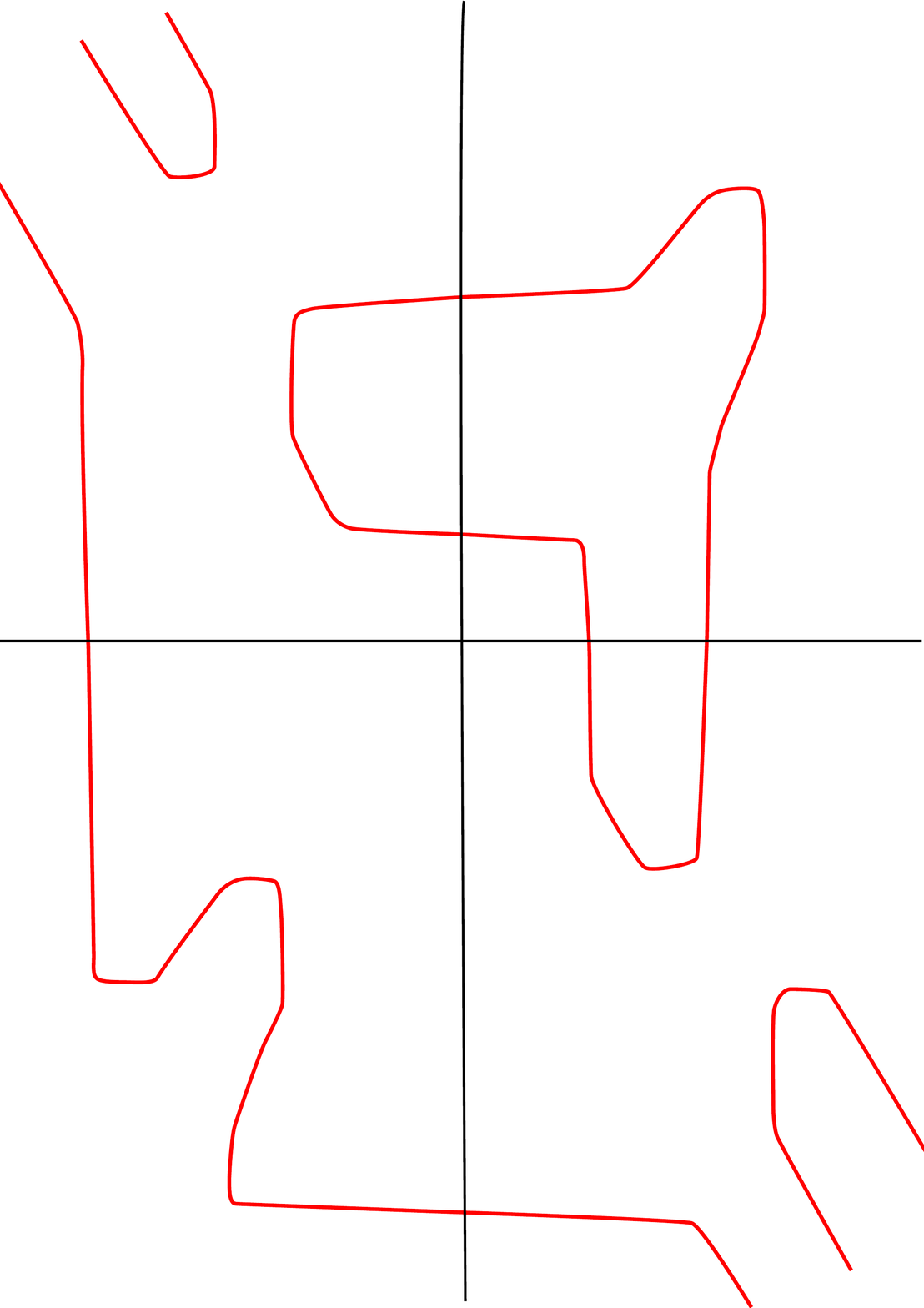}
  \caption{The real part  $\R V$ for the real cubic  from Figure \ref{fig:realcubic} and Example \ref{ex:realcubic}. \label{fig:realcubic2} }
  \end{figure}
\end{example}

Let $\tilde{C}$ denote the first barycentric subdivision of $C$, which results in adding a vertex in the middle of each edge. 
Then the vertices of  $\tilde{C}$  are the vertices of $C$ together with additional vertices $v_e$ for each edge of $C$.  For every 
edge $e$ of $C$ there are now two edges $e'$ and $e''$ of $\tilde{C}$,  moreover $v_e$ is in the boundary of each of these edges.

We can extend any cellular cosheaf $\mathcal{G}$, in particular, $\F_0, \F_1,$ or $\mathcal{S}_{\mathcal{E}}$, to a cellular cosheaf on $\tilde{C}$ in the following way. 
Set $\mathcal{G}(\tilde{e}')= \mathcal{G}(\tilde{e}'') = \mathcal{G}(e)$. 
If $v_e$ is the midpoint of an edge $e$ then  define $\mathcal{G}(v_e) = \mathcal{G}(e)$. The cosheaf morphisms $\mathcal{G}(\tilde{e}') \to \mathcal{G}(v_e) $ are the identity maps. 
Changing the  cellular structure   does not change the homology groups of the cosheaves $\F_0, \F_1,$ and $\mathcal{S}_{\mathcal{E}}$. Namely, 
$H_i(\tilde{C};\F_0) \cong H_i({C};\F_0)$,  $H_i(\tilde{C};\F_1) \cong H_i({C};\F_1),$ and $H_i(\tilde{C} ;\mathcal{S}_{\mathcal{E}}) \cong H_i({C};\mathcal{S}_{\mathcal{E}}).$

For a cellular homology class  $\gamma  \in H_1(C; \F_0)$, we denote by $\text{Supp}(\gamma)$ the collection of edges of $C$ appearing in some chain representing $\gamma$. This is well defined since we are working with $\Z_2$-coefficients.
 
\begin{thm}
\label{thm:twisted}
Let  $C$ be a non-singular compact tropical curve in  a tropical  toric surface. Suppose $C$ is equipped with a real phase structure corresponding to a collection of twists  $T$ of edges of $C$. Then the boundary map of the spectral sequence
 $\partial_1 \colon H_1(\tilde{C};\F_0) \to H_0(\tilde{C};\F_1)$ is given by 
$$
\partial_1(\gamma)=\sum_{e \in T \cap \text{Supp}(\gamma)}  v_e  \otimes s_e,
$$
where $s_e$ is the generator of $\F_1(v_e)$. 
In particular, the number of connected components of $\R C$ is equal to $\dim \Ker (\partial_1) + 1$
\end{thm}
\begin{proof}
It is enough to prove the statement  for cycles in $C$ which are boundaries  of  bounded connected components of the complement $\R^2 \backslash C$ since they form a basis of $H_1(C;\F_0)$. 
Given such a cycle $\gamma \in C_1(\tilde{C};\F_0)$, we first choose a lift $\tilde{\gamma} \in C_1(C;\S_\mathcal{E})$ as follows.  Let  $v$  be a trivalent vertex of $C$ and suppose that $v $ is in the cycle 
$\gamma$. Let $\tilde{e}_1 $ and $\tilde{e}_2$ be the two edges of $\tilde{C}$  (or half edges in $C$) which  share the endpoint $v$ and are contained in $\gamma$, see Figure \ref{fig:Haas}. Let $\varepsilon(v)$ denote the unique element in  $\mathcal{S}_{\mathcal{E}}(\tilde{e}_1)\cap \mathcal{S}_{\mathcal{E}}(\tilde{e}_2)$ by Definition \ref{realphase}.
\begin{figure}
  \centering \includegraphics[width=5cm,height=5cm]{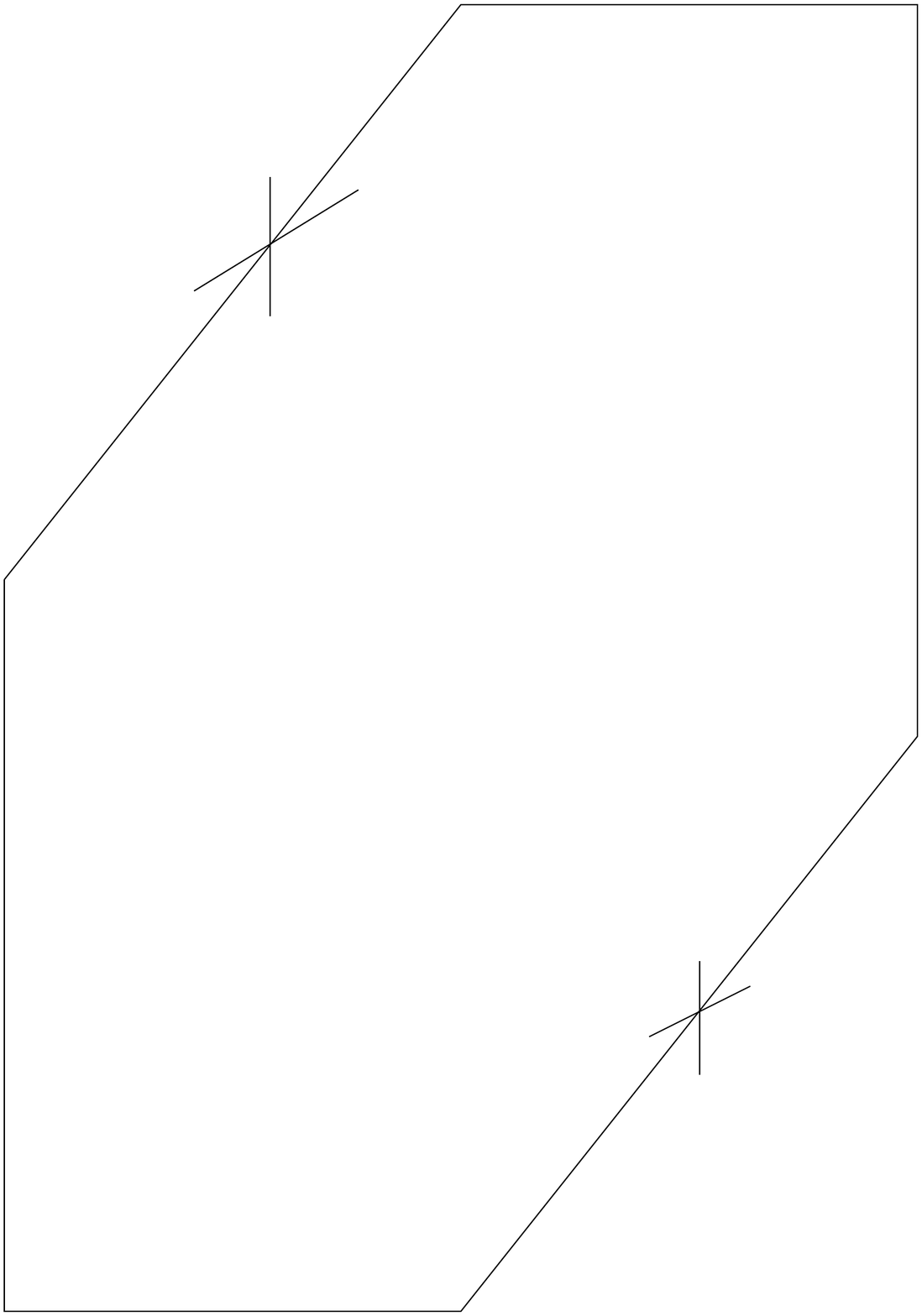}
 \put(-142,78){\tiny $\bullet$} 
 \put(-142,40){\tiny $\bullet$}
 \put(-102,114){\tiny $\bullet$}
 \put(-155,55){$\tilde{e}_1$}
 \put(-145,100){ $\tilde{e}_2$}
 \put(-150,78){$v$}
 \put(-138,60){\tiny $\{\va',\va(v)\}=\mathcal{S}_{\mathcal{E}}(e_1)$}
 \put(-124,90){\tiny $\{\va'',\va(v)\}=\mathcal{S}_{\mathcal{E}}(e_2)$}
 \caption{The cycle $\gamma$ of the cubic from Figure \ref{fig:realcubic} and the lift around a vertex. \label{fig:Haas}} 
\end{figure} 

We set $$\tilde{\gamma} = \sum_{\tilde{e} \in  \gamma \cap \tilde{C}} \tilde{e} \otimes w_{\varepsilon(v)} \in C_1(C; \mathcal{S}_{\mathcal{E}}), $$
where in the sum above $v$ is the unique trivalent vertex of $\tilde{C}$ adjacent to the edge $\tilde{e}$.

If $e\in\mathrm{Edge}(C \cap \gamma)\cap T$ and $v,v'$ are the two adjacent vertices of $e$, then  $w_{\varepsilon(v)}$ and $w_{\varepsilon(v')}$ are different and 
$$
bv_1(w_{\varepsilon(v)}+w_{\varepsilon(v')})=w_e \in \F_1(e).
$$
It $e$ is not twisted, then $w_{\varepsilon(v)}=w_{\varepsilon(v')}$. This proves that $\partial\tilde{\gamma}\in C_0(C;\K_1)$ is supported by the midpoints of twisted edges and that the image by $bv_1$ of the coefficient over $e$ is exactly the generator of $\F_1(e)$. This proves the lemma.
\end{proof}

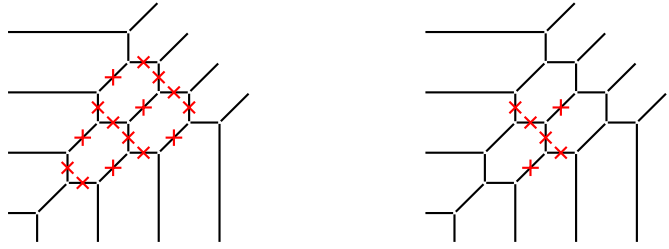
\begin{figure}
 \begin{center}
 
\begin{tikzpicture} [scale = .4, thick = 1mm, every node/.style={inner sep=0,outer sep=0}]


\node (i) at (0,0) {};
\node (i0) at (-1,0)  {} ;
\node (i1) at (0,-1)  {};
\node (i2) at (1,1    )  {};

 \foreach \from/\to in {i/i0,i/i1,i/i2}
  \draw[black] (\from) -- (\to);

\node (j) at (1,2)    {};
\node (j0) at (-1,2)  {} ;
\node (j1) at (2,3)   {};
\node (j2) at (1,1    )   {};

 \foreach \from/\to in {j/j0,j/j1,j/j2}
  \draw[black] (\from) -- (\to);

\node (j) at (2,1)    {};
\node (j0) at (2,-1)  {} ;
\node (j1) at (3,2)   {};
\node (j2) at (1,1    )   {};
 \foreach \from/\to in {j/j0,j/j1,j/j2}
  \draw[black] (\from) -- (\to);

\node (j) at (2,4)    {};
\node (j0) at (-1,4)  {} ;
\node (j1) at (2,3)   {};
\node (j2) at (3,5    )   {};

 \foreach \from/\to in {j/j0,j/j1,j/j2}
  \draw[black] (\from) -- (\to);

\node (j) at (4,2)    {};
\node (j0) at (4,-1)  {} ;
\node (j1) at (3,2)   {};
\node (j2) at (5,3    )   {};

 \foreach \from/\to in {j/j0,j/j1,j/j2}
  \draw[black] (\from) -- (\to);

\node (j) at (3,3)    {};
\node (j0) at (2,3)  {} ;
\node (j1) at (3,2)   {};
\node (j2) at (4,4    )   {};

 \foreach \from/\to in {j/j0,j/j1,j/j2}
  \draw[black] (\from) -- (\to);

\node (j) at (3,6)    {};
\node (j0) at (-1,6)  {} ;
\node (j1) at (4,7)   {};
\node (j2) at (3,5    )   {};

 \foreach \from/\to in {j/j0,j/j1,j/j2}
  \draw[black] (\from) -- (\to);

\node (j) at (6,3)    {};
\node (j0) at (6,-1)  {} ;
\node (j1) at (7,4)   {};
\node (j2) at (5,3    )   {};

 \foreach \from/\to in {j/j0,j/j1,j/j2}
  \draw[black] (\from) -- (\to);

\node (j) at (4,5)    {};
\node (j0) at (3,5)  {} ;
\node (j1) at (4,4)   {};
\node (j2) at (5,6    )   {};

 \foreach \from/\to in {j/j0,j/j1,j/j2}
  \draw[black] (\from) -- (\to);

\node (j) at (5,4)    {};
\node (j0) at (5,3)  {} ;
\node (j1) at (4,4)   {};
\node (j2) at (6,5    )   {};

 \foreach \from/\to in {j/j0,j/j1,j/j2}
  \draw[black] (\from) -- (\to);

  \tikzset{cross/.style={cross out, draw=black, minimum size=2*(#1-\pgflinewidth), inner sep=0pt, outer sep=0pt},
cross/.default={3pt}}
\draw (1,1.5) node[cross, red] {};
\draw (1.5,1) node[cross,red] {};

\draw (2.5,1.5) node[cross,rotate = 45,red] {};
\draw (1.5,2.5) node[cross,rotate = 45,red] {};

\draw (2.5,3) node[cross,red] {};
\draw (3,2.5) node[cross,red] {};
\draw (2,3.5) node[cross,red] {};
\draw (3.5,2) node[cross,red] {};

\draw (2.5,4.5) node[cross,rotate = 45,red] {};
\draw (4.5,2.5) node[cross,rotate = 45,red] {};
\draw (3.5,3.5) node[cross,rotate = 45,red] {};

\draw (3.5,5) node[cross,red] {};
\draw (5,3.5) node[cross,red] {};

\draw (4.5,4) node[cross,red] {};
\draw (4,4.5) node[cross,red] {};

\end{tikzpicture}
\hspace{2cm}
\begin{tikzpicture} [scale = .4, thick = 1mm, every node/.style={inner sep=0,outer sep=0}]


\node (i) at (0,0) {};
\node (i0) at (-1,0)  {} ;
\node (i1) at (0,-1)  {};
\node (i2) at (1,1    )  {};

 \foreach \from/\to in {i/i0,i/i1,i/i2}
  \draw[black] (\from) -- (\to);

\node (j) at (1,2)    {};
\node (j0) at (-1,2)  {} ;
\node (j1) at (2,3)   {};
\node (j2) at (1,1    )   {};

 \foreach \from/\to in {j/j0,j/j1,j/j2}
  \draw[black] (\from) -- (\to);

\node (j) at (2,1)    {};
\node (j0) at (2,-1)  {} ;
\node (j1) at (3,2)   {};
\node (j2) at (1,1    )   {};
 \foreach \from/\to in {j/j0,j/j1,j/j2}
  \draw[black] (\from) -- (\to);

\node (j) at (2,4)    {};
\node (j0) at (-1,4)  {} ;
\node (j1) at (2,3)   {};
\node (j2) at (3,5    )   {};

 \foreach \from/\to in {j/j0,j/j1,j/j2}
  \draw[black] (\from) -- (\to);

\node (j) at (4,2)    {};
\node (j0) at (4,-1)  {} ;
\node (j1) at (3,2)   {};
\node (j2) at (5,3    )   {};

 \foreach \from/\to in {j/j0,j/j1,j/j2}
  \draw[black] (\from) -- (\to);

\node (j) at (3,3)    {};
\node (j0) at (2,3)  {} ;
\node (j1) at (3,2)   {};
\node (j2) at (4,4    )   {};

 \foreach \from/\to in {j/j0,j/j1,j/j2}
  \draw[black] (\from) -- (\to);

\node (j) at (3,6)    {};
\node (j0) at (-1,6)  {} ;
\node (j1) at (4,7)   {};
\node (j2) at (3,5    )   {};

 \foreach \from/\to in {j/j0,j/j1,j/j2}
  \draw[black] (\from) -- (\to);

\node (j) at (6,3)    {};
\node (j0) at (6,-1)  {} ;
\node (j1) at (7,4)   {};
\node (j2) at (5,3    )   {};

 \foreach \from/\to in {j/j0,j/j1,j/j2}
  \draw[black] (\from) -- (\to);

\node (j) at (4,5)    {};
\node (j0) at (3,5)  {} ;
\node (j1) at (4,4)   {};
\node (j2) at (5,6    )   {};

 \foreach \from/\to in {j/j0,j/j1,j/j2}
  \draw[black] (\from) -- (\to);

\node (j) at (5,4)    {};
\node (j0) at (5,3)  {} ;
\node (j1) at (4,4)   {};
\node (j2) at (6,5    )   {};

 \foreach \from/\to in {j/j0,j/j1,j/j2}
  \draw[black] (\from) -- (\to);

  \tikzset{cross/.style={cross out, draw=black, minimum size=2*(#1-\pgflinewidth), inner sep=0pt, outer sep=0pt},
cross/.default={3pt}}

\draw (2.5,3) node[cross,red] {};
\draw (3,2.5) node[cross,red] {};
\draw (2,3.5) node[cross,red] {};

\draw (2.5,1.5) node[cross,rotate = 45,red] {};
\draw (3.5,3.5) node[cross,rotate = 45,red] {};

\draw (3.5,2) node[cross,red] {};
\end{tikzpicture}

\caption{\label{concurrent} The curve $C$ with two collections of twists $T_1$ and $T_2$ from Example \ref{ex:concurrent}. }
\end{center}
\end{figure}

\begin{example}\label{ex:concurrent}
Consider the tropical curve $C$ in both sides of  Figure \ref{concurrent}. The red markings on the edges denote collections of twisted edges $T_1$ on the left  and $T_2$ and the right.It can be verified that both collections of twists are admissible. 

 Consider the basis  $\gamma_1, \gamma_2, \gamma_3 $ of $H_1(C; \F_0)$ where $\gamma_i$'s are the boundaries of the three bounded connected components of $\R^2 \backslash C$. Let $\gamma_1^*, \gamma_2*, \gamma_3^*$ denote the dual basis of $H_0(C; \F_1)$. 
We can represent the map from $\partial_1$ from Theorem \ref{thm:twisted} by a matrix using these two ordered bases and we obtain the matrices 
$$\left(\begin{array}{ccc}0 & 1 & 1 \\1 & 0 & 1 \\1 & 1 & 0\end{array}\right) \qquad \text{and} \qquad \left(\begin{array}{ccc}1 & 1 & 1 \\1 & 1 & 1 \\1 & 1 & 1\end{array}\right),$$
for the twists $T_1$ and $T_2$, respectively. 
The matrix on the left has a $1$-dimensional null space, and therefore a real algebraic curve produced from the collection of twists on the left of Figure \ref{concurrent}  has two connected components. On the right the matrix has a $2$-dimensional null space and the curve from the twists on the right of Figure \ref{concurrent}  has $3$ connected components.  
\end{example}

\subsection{M-curves and Haas theorem}
Haas in his thesis \cite{Haasthesis} studied maximal curves obtained by primtive patchworking. In particular, he found a necessary and sufficient criterion for maximality (see also \cite[Section 3.3]{BIMS} and \cite{BBR}).
Here as an example we reformulate and reprove Haas' criterion for maximality using the techniques of the last section. 

\begin{defi} 
An edge $e$ of a plane tropical curve $C$ is called \emph{exposed} if $e$ is in the closure of an  unbounded connected component of $\R^2 \backslash C$. The set of exposed edges is denoted by $\mathrm{Ex}(C)$. Denote by $\mathrm{Ex}^c(C)$ the  complement of $\mathrm{Ex}(C)$  in the set of bounded edges of $C \cap \R^2$.
\end{defi}

The following theorem is a reformulation of Haas' maximality condition reproved using the description of the map in the spectral sequence from Theorem \ref{thm:twisted}.

\begin{thm}[Haas' maximality condition \cite{Haasthesis} ]\label{HaasMax}

A non-singular compact tropical curve $C$  in a tropical toric variety equipped with a real phase structure corresponding to a 
 collection of twisted edges  $T \subset \text{Edges}(C)$ produces a 
 maximal curve  if and only if $T\cap\mathrm{Ex}^c(C)=\emptyset$ and for every cycle $\gamma \in H_1(C;\Z_2)$ the intersection  
$\gamma \cap T$ consists of an even number of edges. 
\end{thm}
\begin{proof}
By Theorem \ref{thmmax}, the curve $\R C$ is maximal if and only if $\partial_1=0$. 
Cycles in $C_1(C ; \F_0)$ which are boundaries of connected components of the complement $\R^2 \backslash C$ form a basis of $H_1(C; \mathcal{F}_0)$. 
There are $g:= b_1(C)$ such cycles and we denote them by $\gamma_1, \dots, \gamma_g$. Therefore, it suffices to show that $\partial_1(\gamma_i) = 0$ for all $i$. 

For $C$ a non-singular tropical curve there is a non-degenerate pairing: 
$$ \langle  \ , \ \rangle \colon H_0(C; \F_1) \times H_1(C; \F_0) \to \Z_2$$
induced from the pairing on integral homology groups for non-sinuglar tropical curves in \cite{ShawThesis}. A similar non-degenerate pairing defined between tropical homology and cohomology groups is also defined in \cite[Section 7.8]{BIMS} and \cite[Section 3.2]{MZ}.   
On the chain level this pairing is:
$$\langle \beta , \gamma \rangle =  |\text{EdgeSupp}(\beta') \cap \gamma | \mod 2,$$
where $\beta' \sim \beta$ and $\beta' \in C_0(\tilde{C}; \F_1)$ is supported on the midpoints of edges of $C$. The set $\text{EdgeSupp}(\beta')$
consists of the edges of $C$ whose midpoint is in the support of $\beta'$. 
Therefore, it suffices to show that for all pairs of such cycles $\gamma_i$ and $\gamma_j$ the non-degenerate pairing    $\langle \partial_1 (\gamma_i) , \gamma_j \rangle $ is zero. 

The intersection $\gamma_i \cap T$  is even if and only if $\langle \partial_1 (\gamma_i) , \gamma_i \rangle  = 0$. 
Secondly, the pairing  $\langle \partial_1 (\gamma_i) , \gamma_j \rangle  = 0$ if and only if $\gamma_i \cap \gamma_j \cap T$ is a set of even  cardinality. 
Since $\gamma_i$ and $\gamma_j$ are boundaries of convex regions in $\R^2$ they can only intersect in at most one edge of $C$. Therefore, the intersection $\gamma_i \cap \gamma_j \cap T$ must be empty and the statement is proved. 
\end{proof}

\bibliographystyle{alpha}
\bibliography{biblio}
\end{document}